\documentclass[cvs,envcountsect]{svjour}

\usepackage{amsfonts}
\usepackage{amssymb}
\usepackage{times}
\usepackage{bbm}

\usepackage[dvips]{graphicx}

\def\Image{\hbox{\rm Image}}
\def\R{\mathbb R}

\def\intx{\int_{x_{i-1}}^{x_i}}
\def\inty{\int_{\tilde x_{i}}^{\tilde x_{i+1}}}

\newcommand{\gradu}{\ensuremath{\nabla u}}

\newcommand{\gradw}{\ensuremath{\nabla w}}
\newcommand{\Partial}[3]{\ensuremath{\frac{\partial^{#1}{#2}}{\partial{#3}^{#1}}}}
\renewcommand{\Partial}[3][ ]{\ensuremath{\frac{\partial^{#1}{#2}}{\partial{#3}^{#1}}}}
\newcommand{\Div}{\ensuremath{ \hbox{div} }}

\begin{document}

\title{Comparison study for Level set and Direct Lagrangian methods
for computing
Willmore flow of closed planar curves}

\author{Michal Bene\v s \inst{1}, Karol Mikula \inst{2}, Tom\'a\v s Oberhuber
\inst{1} 
\and Daniel \v Sev\v covi\v c \inst{3}
\thanks
{
The authors were  partly supported by the following projects and grants: 
the project HPC-EUROPA(RII3-CT-2003-506079), the NCMM project LC06052, 
VEGA 1/3321/06 grant, the project MSM 6840770010 and APVV-0247-06,
APVV-RPEU-0004-06 grants.
}
}
 
\institute{
Department of Mathematics, Faculty of Nuclear Sciences and Physical Engineering,
Czech Technical University in Prague, Trojanova 13, 120 00 Praha 2, Czech Republic
\vglue0mm
\email{benes@kmlinux.fjfi.cvut.cz},
\email{oberhuber@kmlinux.fjfi.cvut.cz}
\and
Department of Mathematics, Slovak University of
Technology, Rad\-lin\-sk\'eho 11, 813 68 Bratislava, Slovak Republic
\vglue0mm 
\email{mikula@vox.svf.stuba.sk}
\and
Dept. of Applied Mathematics and Statistics,
Faculty of Mathematics, Physics \& Informatics, Comenius University, 842 48
Bratislava, Slovak Republic
\vglue0mm
\email{sevcovic@fmph.uniba.sk}
}

\maketitle
 
\begin{abstract}
The main goal of this paper is to present results of comparison study for the level set and direct Lagrangian methods for computing evolution of the Willmore flow of embedded planar curves.
To perform such a study we construct new numerical approximation schemes for both Lagrangian as well as level set methods based on semi-implicit in time and finite/complementary volume 
in space discretizations. The Lagrangian scheme is stabilized in tangential
direction  by the asymptotically uniform grid point redistribution. Both methods are
experimentally  second order accurate. Moreover, we show precise coincidence of both
approaches in case of various elastic curve evolutions provided that
solving the linear systems in semi-implicit level set method 
is done in a precise way, redistancing is performed occasionally
and the influence of boundary conditions on the level set function is eliminated.
\end{abstract}


\thispagestyle{plain}
\markboth{M.BENE\v S, K. MIKULA, T. OBERHUBER, D.\v SEV\v COVI\v C}
{COMPARISON STUDY FOR WILLMORE FLOW OF PLANAR CURVES}

\section{Introduction}

In the past years, elastic curves, the Willmore functional and the corresponding
gradient flow (the Willmore flow) attracted a lot of attention from both
theoretical as well as computational point of view. Following Daniel Bernoulli's
model of an elastic rod, a classical elastica is
a curve $\Gamma$ in the plane which is a critical point (minimizer) for the
elastic energy functional 
\begin{equation}
E(\Gamma) = \frac12 \int_\Gamma k^2 \,\hbox{d}s\,.
\label{elastic-energy}
\end{equation}
The first comprehensive study of analytical properties of non-closed planar
curves that are minimizers to  (\ref{elastic-energy}) goes back to Leonhard
Euler who presented their characterization and classification in the pioneering
work Additamentum I (De Curvis Elasticae) contained in his Opera Omnia \cite{E}.
Since then much effort has been spent to analyze and provide complete
characterization of both minimizers to (\ref{elastic-energy}) as well as
solutions corresponding to the gradient flow associated with the elastic
energy functional (\ref{elastic-energy}). It is well known from  Euler's work 
that the flow of planar curves with the normal velocity given by
\begin{equation}
\label{beta-eq}
\beta = - \partial^2_s k -\frac12 k^3
\label{geomrov}
\end{equation}
is a gradient flow for the elastic energy functional $E(\Gamma)$ (see e.g. 
\cite{DG,DKS}).
Such fourth order flows of closed curves and its 3D analogies appear in 
various physical 
and computer vision applications dealing with a motion of phase interfaces
or with an image and surface reconstructions \cite{CT,DKS,Se2,KWT,CDDRR,CS,ZC}.

We remind ourselves that the so-called surface diffusion problems (see e.g.
\cite{BMN,MS_ALG}) are described by nonstationary $4^{th}$ order intrinsic
partial differential equations. Similarly, a numerical solution to the Willmore
flow, either in direct (Lagrangian) or level set (Eulerian) formulation, is
a nontrivial problem and leads to a solution of fourth order in space nonlinear
evolution PDEs  Convergence of a semidiscrete time continuous finite element discretization 
in the case when the evolved surface is a graph has been proved by Dziuk and 
Deckelnick in \cite{DD}. First numerical study based on the finite element method for the
Willmore flow in Lagrangian formulation was presented in \cite{DKS} and for the
level set formulation in \cite{DR}. Finite difference discretization has been 
analyzed in \cite{B,O}. Tangential stabilization of Lagrangian
approach for solving fourth order elastic curve flows in case of surface diffusion
was first introduced in \cite{MS_ALG}. Then a parametric finite element method
was tangentially stabilized in \cite{BGN}. Although the Lagrangian methods are
fast and robust (when incorporating a suitable tangential velocity) they cannot
handle topological changes for which the level set methods are
preferred \cite{Se2,DR}. However, a careful and systematic comparison of
nontrivial examples of direct and level set approaches for fourth order curve
evolution problems is still missing. The goal of this paper is to provide such a
comparison study, and, moreover to derive new numerical schemes based on the
finite/complementary volume strategies for both Lagrangian and level
set formulations of the Willmore flow.

First, we present a tangentially stabilized Lagrangian met\-hod based on a
solution to the curvature, local length and position vector equations
accompanied by the asymptotically uniform tangential redistribution of
numerical grid points. We show experimentally that the method is second order
accurate. We apply this method to various examples of evolution of planar
embedded curves. Stabilization by the tangential velocity allows us to use 
reasonable large computational time steps and prevent formation of various
instabilities like merging of evolving curve representing grid points or swallow
tails, which are typical disadvantages of the direct methods.

Then we develop new semi-implicit complementary volume scheme for solving level
set formulation of the Willmore flow. It is again second order accurate. Due
to a finite volume character of discretization it has a potential to be
naturally connected with finite volume schemes for advective level set equations
\cite{FM1,FM2} and thus to be used in various models where the fourth order
terms serve as a curve motion regularization arising, e.g., in image
segmentations \cite{ZC}.
 
The outline of the paper is as follows. In section 2.1 we recall a closed
governing system of equations for the curvature,  local length and position
vector describing evolution of plane curves satisfying (\ref{geomrov}) in
Lagrangian formulation and describe the main idea of asymptotically uniform
tangential redistribution. Section 2.2 focuses on the brief derivation of the
governing equation representing the evolution of zero level set satisfying 
the geometric equation (\ref{geomrov}). In section 3.1 we present our
Lagrangian numerical approximation scheme and, in section 3.2, approximation of
the level set equation for the Willmore flow. Section 4 is devoted to study of
the experimental order of convergence for both methods and to comparison of both
methods in various elastic curve evolution examples.

\section{Governing equations}

\subsection{Direct Lagrangian method}

Henceforth we shall parameterize an embedded regular plane  curve $\Gamma$ 
by a smooth function 
$x:S^1\to \R^2$, i.e. $\Gamma=\Image(x):=\{ x(u), u\in S^1 \}$ such that the
local 
length element $g=|\partial_u x| >0$ is everywhere positive. Taking into account the
periodic boundary conditions at 
$u=0,1$ we shall hereafter identify $S^1$ with the interval $[0,1]$. The unit 
arc-length parameterization will be denoted by 
$s$, so $\hbox{d} s = g\, \hbox{d} u$. 
The tangent vector $\vec T$ and the signed curvature 
$k$ of $\Gamma$ satisfy 
$\vec T = \partial_s x= \partial_u x / g$,
$k = \partial_s x\wedge \partial^2_s x =
\partial_u x\wedge \partial^2_u x / g^3$.
Moreover, we choose the unit inward normal vector $\vec N$  such that 
$\vec T \wedge \vec N =1$ where $\vec a \wedge \vec b$ is the determinant
of the $2\times 2$ matrix with column vectors $\vec a,\vec b$. 
By $\nu$ we denote the tangent angle to $\Gamma$, i.e.
$\vec T = (\cos \nu, \sin \nu )^T$. Now it follows from Fren\'et's formulae 
that 
$\partial_s\vec T = k \vec N$, $\partial_s\vec N = -k \vec T$ and
$\partial_s\nu = k$.  Notice that the curvature $k$ is positive for convex
closed curves in our 
convention of picking of  normal and tangent vector orientation.

Let a regular smooth initial curve $\Gamma_0=\Image(x_0)$ be given. According to
\cite{MS2}, an evolving family of planar curves $\Gamma_t = \Image(x(., t)),
t\in [0,T)$, satisfying 
(\ref{geomrov}) can be represented  by a solution 
to the following system of PDEs
\begin{eqnarray}
&&\partial_t k =\partial^2_s \beta +\alpha \partial_s k +k^2\beta\,,
\label{rovnice1} \\
&&\partial_t g  = -g k\beta  + g \partial_s\alpha\,,\label{rovnice3}    \\
&&\partial_t x  = \beta \vec N + \alpha \vec T\,,   \label{rovnice4}
\end{eqnarray}
subject to initial conditions
$k(.,0)=k_0\,,\ g(.,0)=g_0\,,$ and $x(.,0)=x_0(.),$. We impose 
periodic boundary conditions at $u=0,1$. Having recalled the general form
of governing equations we are able to calculate the time derivative of the
elastic energy functional

\[
2 \frac{\hbox{d}}{\hbox{d}t}E(\Gamma_t) = \frac{\hbox{d}}{\hbox{d}t}\int_0^1 k^2 g \,\hbox{d}u = 
\int_{\Gamma_t} 2 k \partial_t k - k^3 \beta +k^2\partial_s\alpha \,\hbox{d}s.
\]
Since $\int_{\Gamma_t} k\partial^2_s\beta = \int_{\Gamma_t}\beta 
\partial^2_s k$ and $\int_{\Gamma_t} \partial_s(\alpha k^2) = 0$ 
we obtain the following equation
\begin{equation}
\frac{\hbox{d}}{\hbox{d}t}E(\Gamma_t) = \int_{\Gamma_t} ( \partial^2_s k + \frac12 k^3) \beta
\,\hbox{d}s\,.
\label{Ederivative}
\end{equation}
It enables us to conclude that the evolution of 
$\Gamma_t$ with the normal velocity 
$\beta=-\partial^2_s k - \frac12 k^3$ is a 
gradient flow (the Willmore flow) for the Willmore elastic energy functional
$E$. 

Notice that the tangential velocity $\alpha$ 
is a free parameter in (\ref{rovnice1})-(\ref{rovnice4}) and it may depend on 
other quantities like e.g. the curvature, normal velocity and/or local length
element 
in various ways including local or nonlocal  dependences, cf.
\cite{Hou1,K2,MS2,MS3,MS_CVS}. 
In this paper we make use of the so-called asymptotically uniform
tangential redistribution derived in 
\cite{MS3,MS_CVS} which is the most natural for the Willmore flow
since an initial shape is approaching evolution of expanding circles. 
Let us denote $L=L_t$ the total length of a curve $\Gamma_t$.  
It follows from analysis of the tangential velocity made in \cite{MS3,MS_CVS}
that
\[
\frac{g(u,t)}{L_t} \to 1  \quad \hbox{as}\ t\to T_{max} 
\quad \hbox{uniformly w.r. to} \ u\in S^1
\]
provided that the tangential velocity is a solution to a non-local equation
\begin{equation}
\partial_s \alpha = k\beta -\langle k\beta \rangle_\Gamma
+ \left(L/g -1\right) \omega,\  \alpha(0,t)=0\,.
\label{alpha-nonlocal}
\end{equation}
Here $\omega>0$ is a given positive constant and $\langle . \rangle_\Gamma$ is
an averaging operator over a curve $\Gamma$, i.e. $\langle k\beta \rangle_\Gamma
=\frac{1}{L}\int_\Gamma k\beta\,\hbox{d}s$. It is clear that redistribution of
grid points along a curve becomes uniform as $t$  approaches the maximal time of
existence $T_{max}$. In the case of a Willmore flow the time horizon is infinite
(i.e. $T_{max}=+\infty$) for planar Jordan curves and  $T_{max}$ can be finite for some
selfintersecting immersed curves in the plane.
Furthermore, inserting $\alpha$ computed from (\ref{alpha-nonlocal}) into
(\ref{rovnice1})--(\ref{rovnice4}) and making 
use of the identity $\alpha \partial_s k = \partial_s(\alpha k) -
k\partial_s\alpha$ then 
the curvature and local length equations can be rewritten as follows
\begin{eqnarray}
&&\partial_t k =\partial^2_s \beta +\partial_s(\alpha k)  
+ k \langle k\beta \rangle_\Gamma  + \left(1-L/g\right) k \omega \,,
\label{prep-rovnice1} \\
&&\partial_t g  = -g \langle k\beta   \rangle_\Gamma 
+ (L-g)\omega\,.\label{prep-rovnice3}
\end{eqnarray}
In other words, the strong "point-wise" influence of the term $k\beta$
in (\ref{rovnice1}) and 
(\ref{rovnice3}) has been softened by the "averaged" term 
$\langle k\beta \rangle_\Gamma$ in
(\ref{prep-rovnice1}) and (\ref{prep-rovnice3}). 
As a consequence, this important property of asymptotically uniform tangential
velocity
enables us to construct an efficient and 
stable numerical scheme preventing fast local decrease
of local lengths (merging of numerical grid points) as well
as forming various further numerical
instabilities related to high local curvature. 
Since
\begin{eqnarray*}
\partial_s^4 x &=& \partial_s^3 \vec T = \partial_s^2 (k \vec N) 
\partial_s^2 k \vec N + 2 \partial_s k \partial_s \vec N + k\partial_s^2
\vec N \\
&=&\partial_s^2 k \vec N - 2 (\partial_s k) k \vec T - k\partial_s(k\vec T)\\
&=&\partial_s^2 k \vec N -3 k (\partial_s k) \vec T - k^2 \partial_s \vec T\\
&=&\partial_s^2 k \vec N -\frac {3}{2}\partial_s (k^2) \partial_s x - k^2
\partial_s^2 x
\end{eqnarray*}
and $\partial^2_s x = k \vec N$ we have 
\begin{eqnarray*}
(-\partial_s^2 k -\frac12 k^3)\vec N &=& 
-\partial_s^4 x - \frac 32 k^2
\partial_s^2 x - \frac {3}{2}\partial_s (k^2) \partial_s x\\
&=&
-\partial_s^4 x  - 
\frac {3}{2} \partial_s(k^2\partial_s x).
\end{eqnarray*}
Thus the governing system of equations (\ref{rovnice1})--(\ref{rovnice4}) for
the Willmore flow
(\ref{geomrov}) with tangential redistribution can be written as follows:
\begin{eqnarray}
&&\partial_t k  =  - \partial^4_s k -\frac 12 \partial^2_s(k^3)
+ \partial_s(\alpha k) + k(k\beta - \partial_s \alpha),
\label{rov2} \\
&&\partial_t \eta = - k\beta + \partial_s \alpha,\ \ \ \eta=\ln(g),
\label{rov1}\\
&&\partial_t x = -\partial_s^4 x - \frac 32 \partial_s(k^2\partial_s x)
+\alpha\partial_s x
\label{rov4}
\end{eqnarray}
where the tangential velocity $\alpha$ is the unique solution to  equation
(\ref{alpha-nonlocal}).

\subsection{Level set method}

In the level set method the evolving family of planar curves $\Gamma_t, t\ge
0,$ is represented by the zero level set of the so-called shape
function $u:\Omega \times [0,T] \to \R$ where $\Omega \subset \R^2$ is a simply
connected domain containing the whole family of evolving curves
$\Gamma_t, t\in[0,T]$. Assuming  zero is the regular value of the mapping $u(.,t)$, 
i.e. $|\nabla u(x,t)|\not=0$ for $u(x,t)=0$ we can express the unit inward normal vector and signed
curvature as: $\vec N =\nabla u/|\nabla u|$ and $k=-\hbox{div
}(\nabla u/|\nabla u|)$. Let us denote the following auxiliary functions:
\[
H=\hbox{div}\left(\frac{\nabla u}{|\nabla u|}\right), \ 
Q=|\nabla u|,\ w=Q H\,.
\]
Then  $\partial_s k = -\nabla H . \partial_s x = -\nabla H
. \vec T$ and, by Fren\'et's formula,  $\partial^2_s k =-k\nabla H . \vec N -
\vec T^t  \nabla^2 H \vec T$. Differentiating the equation $u(x(s,t), t) =0$
with respect to time we obtain $\partial_t u + \nabla u . \partial_t x  =0$.
Since the normal velocity of $x$ is $\beta=\partial_t x . \vec N$ we obtain
$\frac{1}{|\nabla u|}\partial_t u = -\beta$. Inserting expressions
for $\partial^2_s k$ and $k=-H$ we obtain
\[
\frac{1}{Q}\partial_t u = \frac12 \hbox{div}\left(\frac{H^2}{Q}\nabla
u\right) - H^3 
- \vec T^t  \nabla^2 H \vec T\,.
\]
Here $\vec T=(-n_2,n_1)$ where $\vec N=(n_1,n_2)$, i.e. $\vec T$ is the vector $\vec N$
rotated by $-\pi/2$. 
Straightforward calculations show that the right hand side of the above
equation can be rewritten in 
the divergent form. The resulting system of two equations governing the
evolution of the shape function has been derived by Droske and Rumpf in \cite{DR} and 
it reads as follows::
\begin{eqnarray}
\label{levelset}
\partial_t u &=& - Q\, \Div \left( \mathbbm{E}  \gradw - \frac{1}{2}
\frac{w^2}{Q^3} \gradu \right), \\
\label{levelset2}
 w &=& Q\, \hbox{div}\left(\frac{\nabla u}{|\nabla u|}\right)
\end{eqnarray}
where the $2\times 2$ matrix $\mathbbm{E} = \frac{1}{Q}\left( \mathbbm{I} -
\frac{\nabla u}{Q} \otimes \frac{\nabla u}{Q} \right)$ is a projection into a
tangential space of the curve representing the zero level set of $u$.
System of equations (\ref{levelset}--\ref{levelset2}) is subject to the initial
condition
\[
u(x,0)=u^0(x)\,,\quad x\in \Omega
\]
and clamped boundary conditions at $\partial\Omega$, i.e. 
$u(x,t)=0$, $\partial_\nu u(x,t)=0$, $x\in\partial\Omega$. 
The initial function $u^0$ is a signed distance function, i.e.
$u^0(x)=\hbox{dist}(x,\Gamma^0), x\in\Omega$.

\section{Numerical approximation schemes}

\subsection{Numerical approximation of the Lagrangian method}

Our numerical approximation of an evolved curve is represented by discrete
plane points $x_i^j$ where the index $i=1,...,n,$ denotes space 
discretization and the index $j=0,...,m,$ stands for a discrete time
stepping. Due to periodic boundary conditions we use additional values
$x_{-1}^j=x_{n-1}^j$, $x_0^j=x_n^j$, $x_{n+1}^j=x_1^j$, $x_{n+2}^j=x_2^j$. If we
take a uniform division of the time interval $[0,T]$ with a time step
$\tau=\frac{T}{m}$ and a uniform division of the fixed
parameterization interval $[0,1]$ with a step $h=1/n$, a point
$x_i^j$ corresponds to $x(ih,j\tau)$. The systems of difference
equations corresponding to (\ref{alpha-nonlocal}), (\ref{rov2}) -- (\ref{rov4}) 
will be solved for discrete quantities $\alpha_i^j$, $\eta_i^j$, $r_i^j$,
$k_i^j$, $x_i^j$, $i=1,...,n,\ j=1,...,m,$ representing approximations of the
unknowns $\alpha$, $\eta$, $gh$, $k$ and $x,$ respectively. Here $\alpha_i^j$
represents the tangential velocity of a flowing node $x_i^j$, and $\eta_i^j$,
$r_i^j\approx |x_i^j-x_{i-1}^j|$ and $k_i^j$ represent piecewise constant
approximations of the corresponding quantities in the so-called flowing
finite volume $\left[x_{i-1}^j,x_i^j\right]$. In order to derive new position
$x_i^j$ we use corresponding flowing dual volumes $\left[\tilde x_{i-1}^j,\tilde
x_i^j\right]$ where $\tilde x_i^j= \frac{x_{i-1}^j+x_{i}^j}{2}$ with
approximate lengths $q_i^j\approx |\tilde x_i^j-\tilde x_{i-1}^j|$. Our
computational method is simple and natural. At the $j$-th discrete time step, we
first find values of the  tangential velocity $\alpha_i^j$ by discretization of (\ref{alpha-nonlocal}). 
Then the values of $\eta_i^j$ are computed and used 
for updating local lengths $r_i^j$ by discretizing equations (\ref{rov1}). 
Using computed local lengths, the intrinsic derivatives are approximated in (\ref{rov2}),
and (\ref{rov4}), and pentadiagonal systems with 
periodic boundary conditions are constructed and solved for new discrete 
curvatures $k_i^j$ and position vectors $x_i^j$. 

In order to discretize (\ref{alpha-nonlocal}) we integrate it
over flowing finite volume $\left[x_{i-1},x_i\right]$ to obtain
\[
\intx\partial_s\alpha \,\hbox{d}s = \intx
k\beta-\langle k\beta \rangle_\Gamma + \omega\left(L/g - 1\right)
\,\hbox{d}s\,.
\]
Hence
\[
\alpha_i -\alpha_{i-1} = r_i (k_i 
\beta_i- \langle k\beta \rangle_\Gamma) +\omega\left(L/n - r_i\right)
\]
where $\alpha_0=0$. Taking discrete time stepping in the previous relation we
obtain following expression for discrete values of the tangential velocity:
\[
\alpha_i^j = \alpha_{i-1}^j + 
r_i^{j-1}(k_i^{j-1}\beta_i^{j-1} -
B^{j-1})+\omega(\frac{L^{j-1}}{n} - r_i^{j-1})
\]
where, for $i=1, ..., n,$
\[
\beta_i^j = - \frac{1}{r_i^j} \left(\frac{k_{i+1}^j-k_i^j}{q_i^j}-
\frac{k_{i}^j-k_{i-1}^j}{q_{i-1}^j}\right)-\frac{1}{2}
\left(k_i^j\right)^3,
\]
\[
q_i^j =\frac12 \left(r_{i}^j+r_{i+1}^j\right),\ \ 
L^j=\sum\limits_{l=1}^n r_l^j, \ \  B^j=\frac{1}{L^j}
\sum\limits_{l=1}^n  r_l^jk_l^j \beta_l^j,
\]
and $\alpha_0^j=0$, i.e. the point $x_0^j$ is moved in the normal direction.

Now, a similar approximation methodology is applied for equation (\ref{rov1}). Thus 
\[
r_i^{j-1}\frac{\eta_i^j-\eta_i^{j-1}}{\tau} =
- r_i^{j-1}k_i^{j-1}\beta_i^{j-1}+\alpha_i^j-\alpha_{i-1}^j
\]
for $i=1,...,n$. It leads to the update formula for local lengths:
\[
r_i^j=\exp(\eta_i^j),\ \ i=1,...,n,
\]
subject to periodic boundary conditions $r_{-1}^j=r_{n-1}^j,\ \ r_0^j=r_n^j,\
\ r_{n+1}^j=r_1^j, \ \ r_{n+2}^j=r_2^j$. New local lengths are used for
approximation of intrinsic derivatives in the curvature equation
(\ref{rov2}). We obtain
$ \intx\partial_t k \,\,\hbox{d}s = \intx -\partial_s^4 k -\frac
12\partial_s^2(k^3)+\partial_s(\alpha k)\,\,\hbox{d}s + \intx k (k\beta -
\partial_s\alpha)\,\,\hbox{d}s$. Hence
\begin{eqnarray}
r_i\frac{\hbox{d} k_i}{\hbox{d}t}&=& -\left[\partial_s^3 k\right]_{x_{i-1}}^{x_i}-\frac 12
\left[\partial_s(k^3)\right]_{x_{i-1}}^{x_i}  \nonumber \\
&& + \left[\alpha k\right]_{x_{i-1}}^{x_i} +  k_i (r_i k_i\beta_i- (\alpha_i-\alpha_{i-1}))
\label{curvdiscr}
\end{eqnarray}
and taking semi-implicit time stepping, i.e. replacing time derivative by
backward difference and treating linear terms at
the current time level $j$ while the nonlinear terms at the level $j-1$, 
and approximating derivative terms on the boundaries of flowing
finite volumes by finite differences we obtain following
pentadiagonal system with periodic boundary conditions for  new 
discrete values of the curvature:
\begin{eqnarray}
\nonumber\\
a_{i}^j k_{i-2}^j+b_{i}^j k_{i-1}^j+c_{i}^j k_{i}^j+d_i^j k_{i+1}^j
+e_i^j k_{i+2}^j = f_i^j
\label{curvsystem}\\
\nonumber
\end{eqnarray}
for $i=1,...,n$, subject to periodic boundary conditions 
$k_{-1}^j=k_{n-1}^j, k_0^j=k_n^j, k_{n+1}^j=k_1^j, k_{n+2}^j=k_2^j$. 
For completeness, a detailed description of the system coefficients is given in Appendix.

Finally we discretize equation (\ref{rov4}) by integrating
in a dual volume
$\left[\tilde x_{i-1},\tilde x_i\right]$ 
to get
\begin{eqnarray}
\inty\partial_t x\,\hbox{d}s &=& \inty 
-\partial_s^4 x - \frac 32 \partial_s(k^2\partial_s x)
+\alpha\partial_s x
\,\,\hbox{d}s\,,
\nonumber\\
q_i\frac{\hbox{d} x_i}{\hbox{d}t}&=&
\left[-\partial^3_s x - \frac 32 k^2\partial_s x 
\right]_{\tilde x_i}^{\tilde x_{i+1}}
+\alpha_i (\tilde x_{i+1}-\tilde x_i)\,.
\nonumber
\end{eqnarray}
Now replacing the time derivative by the backward difference, derivative terms on boundaries
of dual volume by finite differences and 
$\tilde x_{i}$ by the average of grid points in the last term, 
we obtain two tridiagonal systems  for updating 
the position vector:
\begin{eqnarray}
\nonumber\\
{\cal A}_{i}^j x_{i-2}^j+{\cal B}_{i}^j x_{i-1}^j+
{\cal C}_{i}^j x_{i}^j+ {\cal D}_{i}^j x_{i+1}^j+
{\cal E}_{i}^j x_{i+2}^j={\cal F}_i^j\label{systempos}\\
\nonumber
\end{eqnarray}
$i=1,...,n,$
subject to periodic boundary conditions
$x_{-1}^j=x_{n-1}^j, x_0^j=x_n^j, x_{n+1}^j=x_1^j, x_{n+2}^j=x_2^j$.
The exact form of coefficients $\cal A, B, C, D, E$ can be found in Appendix.

The initial quantities for the algorithm are computed from discrete
representation of the initial curve $x_0$. The reader is referred to
\cite{MS_CVS} for further details. Every pentadiagonal system is solved by
Gauss-Seidel iterates. We stop the Gauss-Seidel iteration procedure if
a difference of subsequent iterates in maximum norm is less than
the prescribed tolerance $10^{-10}$.

\subsection{Numerical approximation of the level set method}

Concerning approximation of the level set equation (\ref{levelset}) we consider
rectangular domain $\Omega \equiv \langle a_1, a_2 \rangle \times \langle b_1,
b_2 \rangle$ and we assume an equidistant spatial step $h$ in both
directions. We define a regular mesh $\omega_h$ consisting of grid points $x_{ij} =
[a_1 + ih, b_1 + jh]$ for $i = 0, ..., N_1, j = 0, ..., N_2, $ where $a_1 +N
_1 h=a_2$ and $b_1 + N_2 h=b_2$. Without loss of generality we shall assume
$a_1 = b_1 = 0$. The corresponding dual mesh $\mathcal V$ is given as the 
union of the finite volumes $V_{ij}$ of the form $\left \langle
\left(i-\frac{1}{2}\right) h, \left(i+\frac{1}{2}\right) h \right\rangle
\times \left \langle \left(j-\frac{1}{2}\right) h,\left(j+\frac{1}{2}\right)
h \right\rangle$ for $i = 0, ..., N_1, j = 0, ..., N_2$. The projection of a
solution at $x_{ij}$ is defined as $u_{ij}=u(x_{ij})$. Similarly as in the
Lagrangian method we take a uniform division of the time interval $[0,T]$ with
a time step $\tau=\frac{T}{m}$. Let us consider an element $V_{ij}$ of the dual
mesh $\mathcal V$. Integrating (\ref{levelset})-(\ref{levelset2}) over $V_{ij}$
and applying the Stokes theorem we obtain
\begin{eqnarray}
\label{willmore-ut-stokes}
\int_{V_{ij}} \frac{1}{Q} \Partial{u}{t} \hbox{d}x&=&
\int_{\partial V_{ij}} \frac{1}{2} \frac{ w^2}{Q^3} \Partial{u}{\nu}
- \left \langle \mathbbm{E} \gradw,  \nu \right\rangle 
\hbox{d}\sigma\,, \\
\label{willmore-w-stokes}
\int_{V_{ij}} \frac{w}{Q}\hbox{d}x &=& \int_{\partial V_{ij}} \frac{1}{Q}
\Partial{u}{\nu} \hbox{d}\sigma
\end{eqnarray}
where $\nu$ is the outer normal of the boundary $\partial V_{ij}$.

We start with approximation of the term $Q$ on $V_{ij}$.
For $r,s \in \{ -1, 1\}, |r|+|s|=1,$ we define the linear operator
$\nabla^{rs}$ as follows:
\begin{eqnarray*}
\nabla^{r,0} u_{ij} &=& \frac{1}{h}\left( r ( u_{i+r,j} -
u_{ij} ), u^{r,1}_{ij} - u^{r,-1}_{ij}  \right), \\ 
\nabla^{0,s} u_{ij} &=& \frac{1}{h} \left( u^{1,s}_{ij} -
u^{-1,s}_{ij} , s (u_{i,j+s} - u_{ij} ) \right)
\end{eqnarray*}
where $u^{rs}_{ij}$ is the average of $u_{ij}$ defined as:
\[
u_{ij}^{rs} = \frac{1}{4} ( u_{ij} + u_{i+r,j} + u_{i,j+s}
+ u_{i+r,j+s} )\,.
\]
For a fixed regularization parameter $0<\epsilon\ll 1$ we define 
\[
Q^{rs; n}_{ij} = \sqrt{ \epsilon^2 + |\nabla^{rs} u^n_{ij}|^2}\,,
\quad \bar Q^n_{ij} = \frac{1}{4} \sum_{|r|+ |s| = 1} Q^{rs; n}_{ij}\,.
\]
\noindent
Let  $\mathbbm{E}_{ij}^{rs; n} =\left( E_{kl; ij}^{rs; n}\right)_{k,l=1,2}$ be
the $2\times 2$ projection matrix: 
\[
\mathbbm{E}^{rs; n}_{ij} = \frac{1}{Q^{rs; n}_{ij}}
\left( \mathbbm{I} - 
   \frac{\nabla^{rs}u^n_{ij}}{Q^{rs; n}_{ij}} \otimes
   \frac{\nabla^{rs}u^n_{ij}}{Q^{rs; n}_{ij}} \right)\,.
\]
Now we are able to derive a discretization of (\ref{willmore-ut-stokes}) 
\begin{eqnarray}
\label{discrete-willmore-u-t}
\frac{u^n_{ij}-u^{n-1}_{ij}}{\tau} &=& 
\frac{\bar Q_{ij}^{n-1}}{2 h^2}\hbox{\hglue-3mm}\sum_{|r|+|s|=1} 
\frac{(\hat w^{rs;n-1}_{ij})^2}{(Q_{ij}^{rs; n-1})^3} 
( u^n_{i+r, j+s} - u^n_{ij}) \nonumber \\
&& -\frac{\bar Q_{ij}^{n-1}}{h^2}\hbox{\hglue-3mm}
\sum_{|r|+|s|=1}\hbox{\hglue-3mm} h\langle \mathbbm{E}_{ij}^{rs, n-1} \nabla^{rs}w^n_{ij},
\nu_{rs} \rangle 
\end{eqnarray}
where
\begin{equation}
\hat w_{ij}^{rs;n} = \frac 12\left(w_{ij}^n+w_{i+r,j+s}^n\right),
\label{w-eq-11}
\end{equation}
and $\nu_{rs}$ is  the unit outer normal vector, $\nu_{rs} = (r,s)$ for $|r|+|s|=1$.
In order to approximate $w^n$ and $w^{n-1}$ on $V_{ij}$ and on its boundary
$\partial V_{ij}$ we have used expression (\ref{willmore-w-stokes}) to obtain
\begin{equation}
w_{ij}^n = \frac{\bar Q_{ij}^n}{h^2} \sum_{|r|+|s|=1}
\frac{1}{Q_{ij}^{rs; n}} (u^n_{i+r,j+s} - u^n_{ij}) \,.
\label{w-eq-1}
\end{equation}
Since (\ref{discrete-willmore-u-t}) contains a new time level $w^n$ expressed
through new time level of the  solution $u^n$ (see  (\ref{w-eq-1})) as well as the previous time level
$w^{n-1}$, the resulting discrete level set scheme is semi-implicit in time. After some
calculations it can be written as twenty one points scheme of the form
\begin{equation}
\label{systemlevelset}
\sum_{(r,s)\in \mathcal{O}} A_{ij}^{rs} u_{i+r,j+s}^n = u_{ij}^{n-1}
\end{equation}
where $\mathcal{O}=\{ (r,s), -2\le r,s\le 2, |r|+|s| < 4\}$ and
 $i = 2, ..., N_1-2, j = 2, ..., N_2-2$. For the remaining $i,j$, 
the values of $u_{ij}^n$ are linearly extrapolated.
The coefficients of the above system can be found in Appendix. 

System (\ref{systemlevelset}) is solved by the iterative GMRES 
algorithm with ILUT (ILU with threshold) preconditioning 
or by the complete LU decomposition as a direct solver (c.f. \cite{SY}). 
The time step $\tau$ is chosen to be proportional to
$h^2$ and the regularization parameter $\epsilon$ can chosen either
as a function of $h$ or $\epsilon$ can be prescribed as a small fixed constant.
As an initial condition we choose a signed distance function $d(x)$ to
the initial curve. At  prescribed redistancing time steps we perform
redistancing of the level set solution back to the signed distance using the fast sweeping method
(see \cite{ZH} for details).
As an alternative to the semi-implicit scheme 
(\ref{discrete-willmore-u-t}) we may also consider its explicit 
version, i.e. all the terms on the right hand side of
(\ref{discrete-willmore-u-t}) are considered at the time step $n-1$,
c.f. \cite{B,O} for other similar explicit schemes. 
In this case we avoid the singularities of the signed distance 
function in its local extrema and the initial condition 
has a "phase-field" like shape 
$$
u^0(x) = \delta \hbox{sgn}(d(x)) ( 1 - \exp(- | d(x) / \delta| ),
$$
where $\delta$ is a parameter describing the width of the region 
where $u^0$ changes from $-\delta$ to $+\delta$.

\section{Discussion on numerical experiments}

\subsection{Experimental order of convergence for the methods}

Let an initial curve be a circle with radius $r_0$. Since for the circle we have
$k = \frac{1}{r}$ then it follows 
from (\ref{beta-eq}) that  $\dot{r}(t) = \frac{1}{2} r(t)^{-3}$. Hence $r(t) =
(2t + r_0^4)^{\frac{1}{4}}$. Using this simple analytical solution we can compute
experimental order of convergence for both schemes. Without loss of generality we choose $r_0=1$. 

In the case of the Lagrangian scheme
we approximate the initial unit circle subsequently by $n=10,20,40$ and $80$
nodes with $h=1/n$. The final time was set up to be $T=2.56$ and time step was chosen to be proportional to $h^2$, i.e. $\tau = h^2$. Table \ref{eoc-tab1} shows errors and experimental 
order of convergence (EOC) of the scheme in $L^p = L^p((0,T),L_p(S^1)) = L^p(S^1\times (0,T) )$ for $p=2, \infty$. In the level set approximation we solve the problem in domain 
$\Omega = \langle -2, 2 \rangle \times
\langle -2, 2 \rangle$. The domain $\Omega$ was  splitted subsequently 
into $n \times n$ finite volumes for $n=10,20,40,80$ with
$h=1/n$. The final time was chosen as $T = 0.5$ and again $\tau = h^2$. The regularization
parameter $\epsilon$ was refined proportionally to the grid refinement
using the rule $\epsilon^2 = 2h$. Finally, the redistancing period was 
$\tau_{redist} = 0.25 h$.  Errors in $L^p, p=2,\infty$ norms are presented in Table \ref{eoc-tab2}.

\begin{table}
\caption{EOC for the Lagrangian scheme in $L^p, p=2,\infty$ norms. }
\begin{tabular}{llllll}
\hline\hline
$\frac{\hbox{Error}}{\hbox{EOC}}$ $\backslash \  h$     & 0.1      & 0.05       &  0.025      & 0.0125   \\ 
\hline\hline
\\
$p=2$ & 0.04301  & 0.01089  & 0.00271 & 0.00067\\
  EOC &          & 1.982  & 2.005 & 2.003 \\
\\
$p=\infty$
      & 0.03402 & 0.00886  & 0.00223  & 0.00056 \\
   EOC&         & 1.940  & 1.986  & 1.988\\
\hline                                                                   
\end{tabular}                                                            
\label{eoc-tab1}
\end{table}

\begin{table}
\caption{EOC for the level set scheme in $L^p, p=2,\infty$ norms.}
\begin{tabular}{llllll}
\hline\hline
$\frac{\hbox{Error}}{\hbox{EOC}}$ $\backslash \  h$     & 0.4       & 0.2        & 0.1        & 0.05   \\ 
\hline\hline
\\
$p=2$ & 0.21497 & 0.06585   & 0.01699  &  0.00400  \\
EOC   &         & 1.707   & 1.954   & 2.086  \\ 
\\
$p=\infty$
      & 0.71190 & 0.12286   & 0.03780   & 0.00973  \\
EOC   &         & 2.534   & 1.700   & 1.957\\
\hline
\end{tabular}
\label{eoc-tab2}
\end{table}

\subsection{Comparison of the Lagrangian and the level set evolutions}

In this section we compare the numerical results obtained by our Lagrangian
and the level set approaches on various representative examples. 
In the case of Lagrangian scheme we approximate an evolving curve
by $100$ grid nodes in all experiments to follow. In the case of  the level 
set method we hereafter split domain $\Omega$ into $100 \times 100$ finite volumes.

\begin{figure}
\center{
\includegraphics[width=3cm]{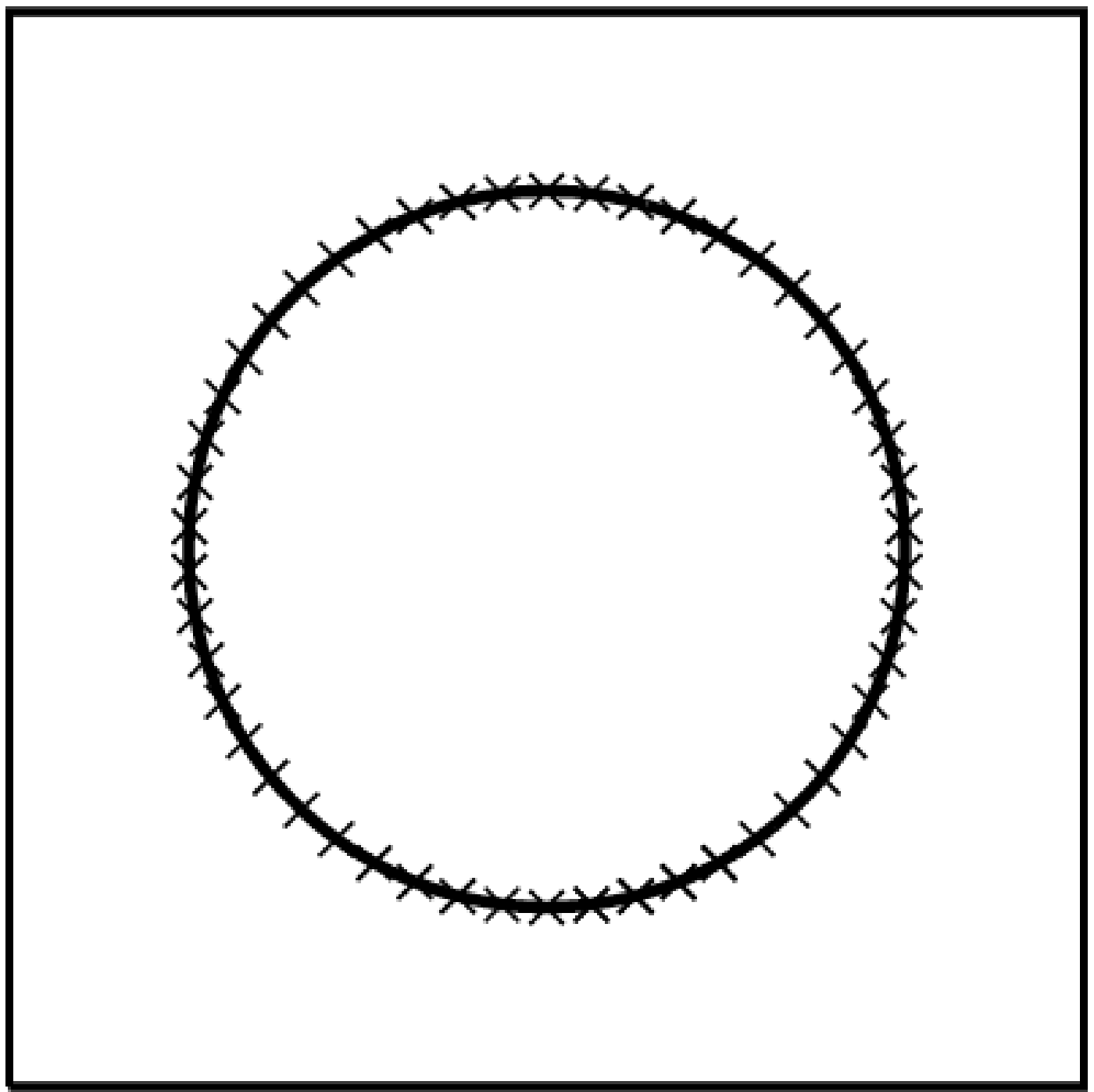}
\hglue 0.5cm
\includegraphics[width=3cm]{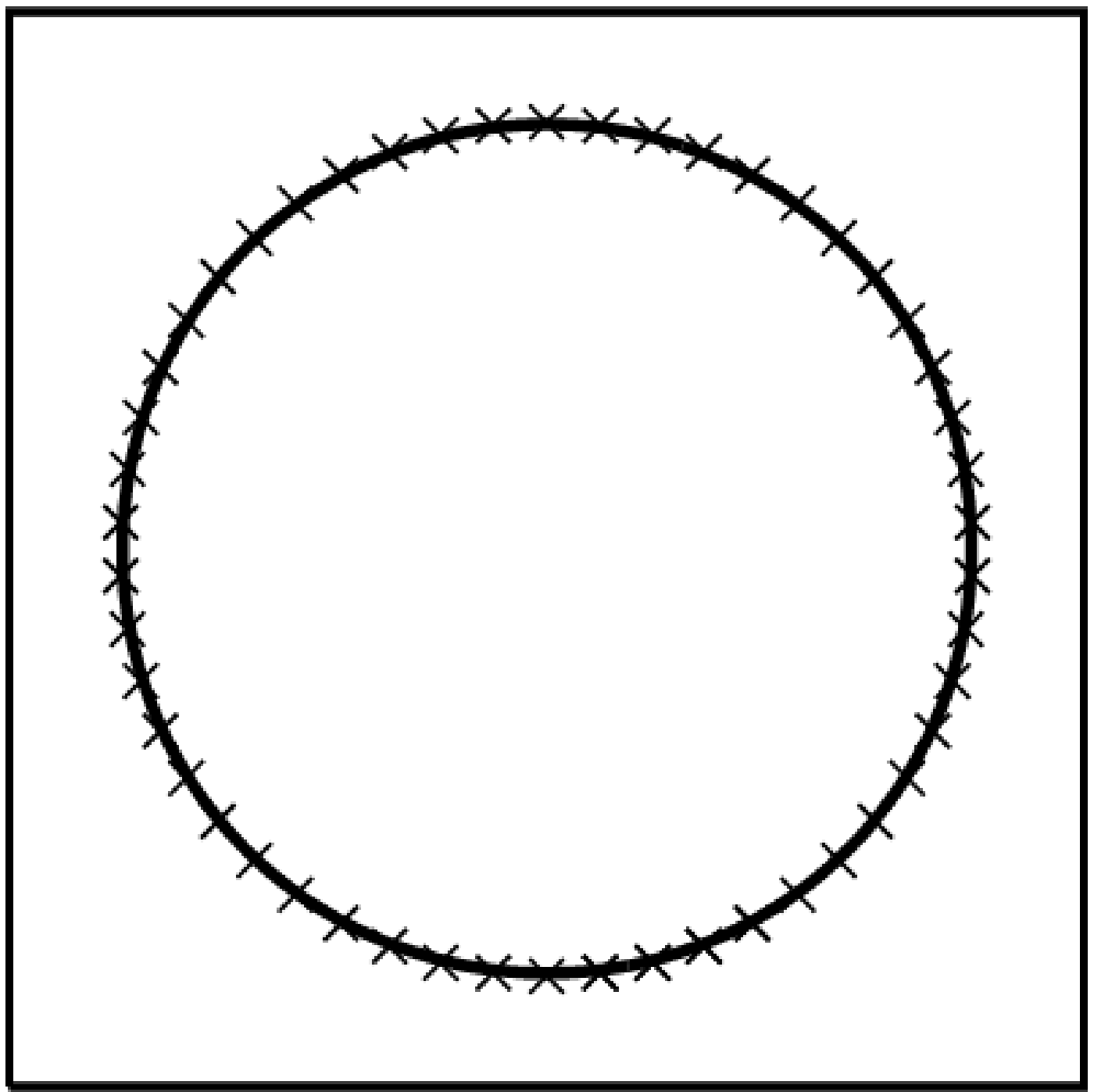}
}
\center{a) \hglue 3truecm b)}

\center{
\includegraphics[width=5.2cm]{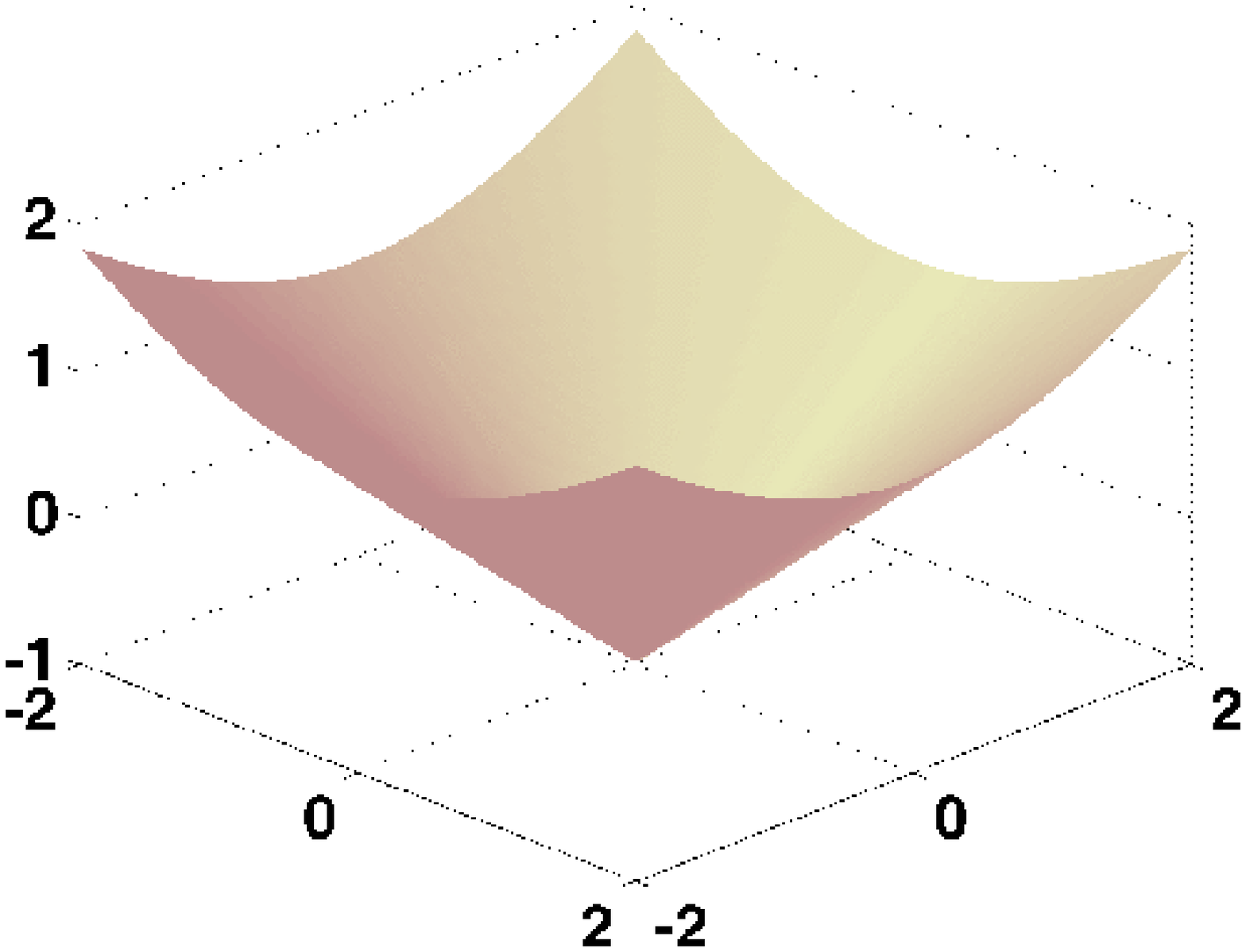}
}
\center{c)}
\center{
\includegraphics[width=5.2cm]{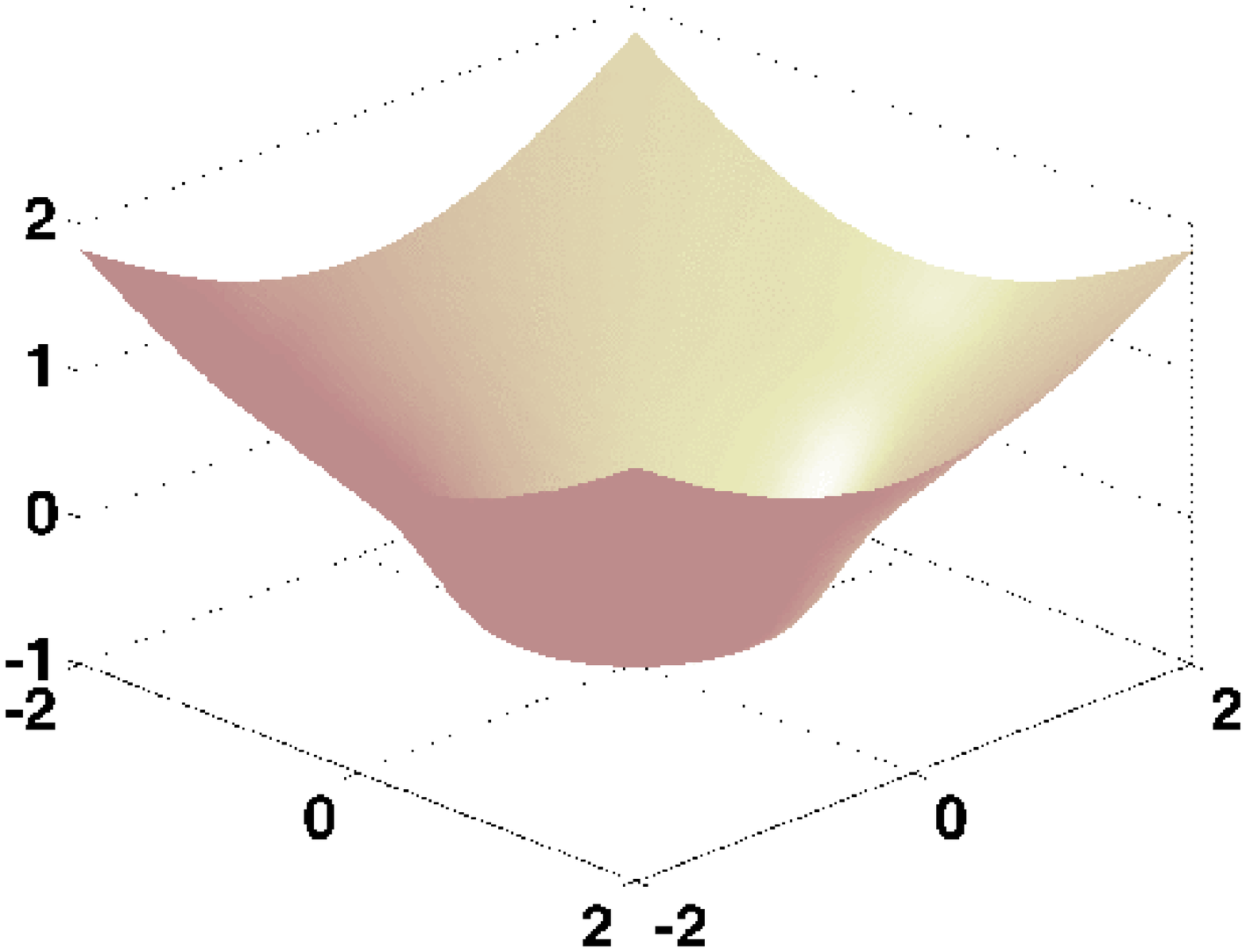}
}
\center{d)}
\caption{
a) A circle as an initial condition;
c) the same initial condition for the level set approximation;
b) a circle computed at the time $T=0.5$ by both methods
(cross marks correspond to  the Lagrangian method);
d) the level set function at the time $T=0.5$.}
\label{circle-test}
\end{figure}

In Fig. \ref{circle-test} we present comparison of both methods for the case of  evolution 
of an initial circle with the radius $r_0=1$. The time horizon $T=0.5$. By cross marks we depict 
approximation by the Lagrangian direct scheme where the evolution was computed using the time step $\tau=0.002$ and no tangential redistribution ($\alpha=0$).
The evolution of the level set function was computed in the spatial domain 
$\Omega = \langle -2, 2 \rangle \times
\langle -2, 2 \rangle$ with the  time step
$\tau=0.002$ and the smoothing parameter  $\epsilon=0.001$.
We did not provide redistancing in this case in order to show 
deformation of an initial distance function to final shape of the level set 
function (see Fig. \ref{circle-test} bottom).

\begin{figure}
\center{
\includegraphics[width=3cm]{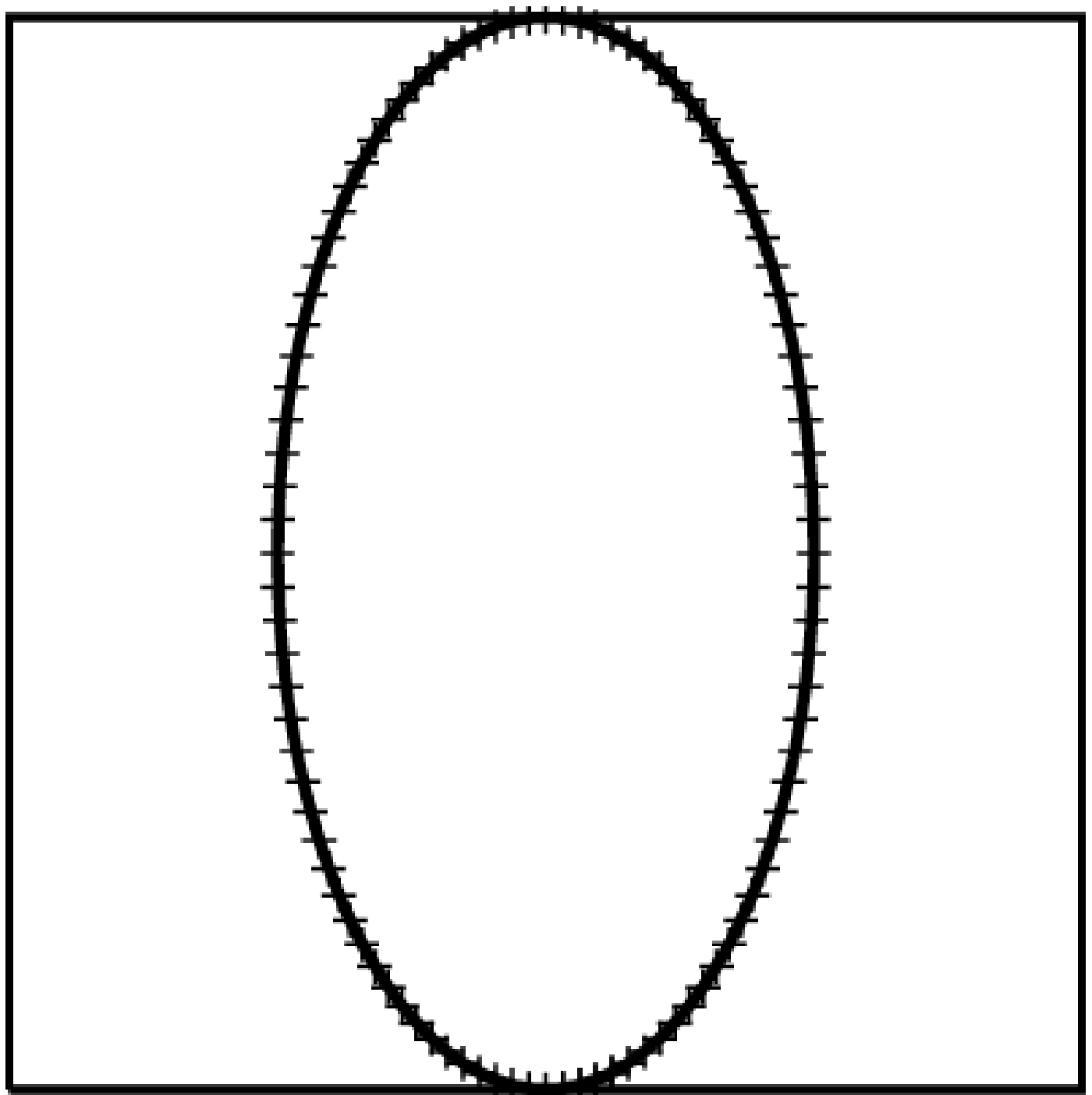}
\hglue 0.5cm
\includegraphics[width=3cm]{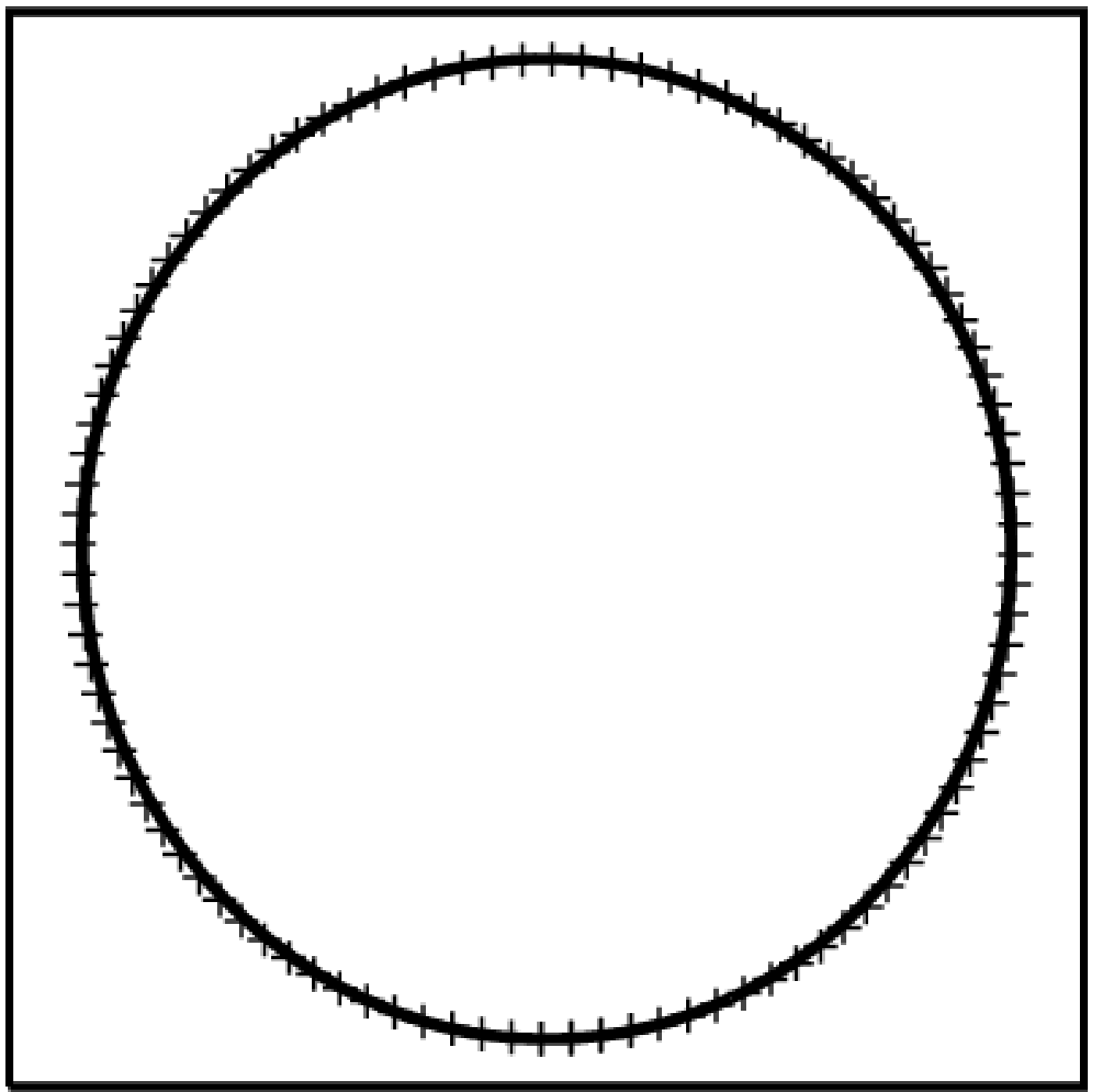}
}
\center{a) \hglue 3truecm b)}

\center{
\includegraphics[width=5.2cm]{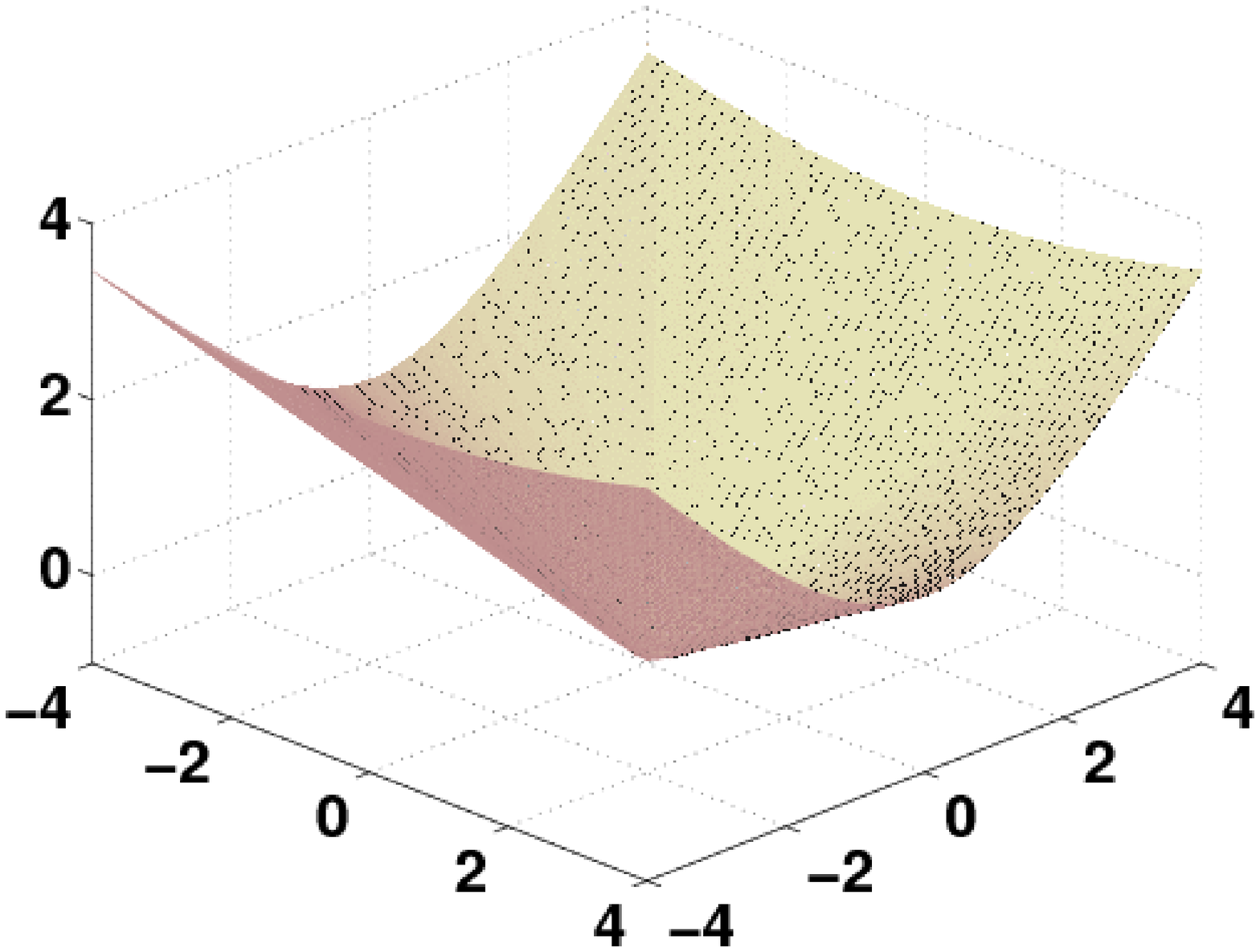}
}
\center{c)}
\center{
\includegraphics[width=5.2cm]{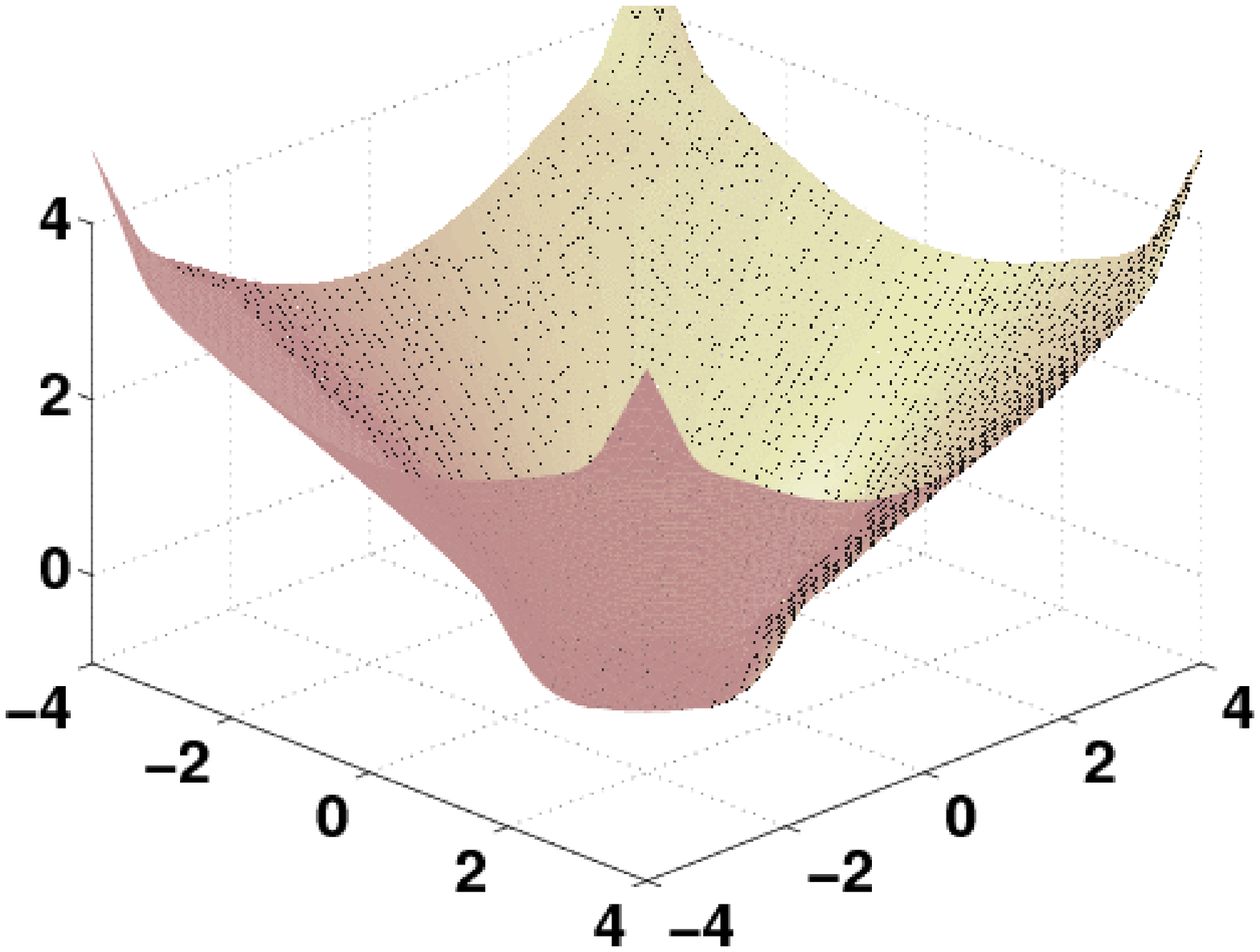}
}
\center{d)}

\caption{
a) An initial ellipse with half-axes ratio 1:2;
c) the same initial condition for the level set approximation;
b) the approximate circle computed at the time $T=2.56$ by both methods
(cross marks correspond to  the Lagrangian method);
d). the level set function at the time $T=2.56$.}
\label{ellipse-test}
\end{figure}

In Fig. \ref{ellipse-test} we show evolution of an initial ellipse with half-axes $1$ and $2$. It asymptotically approaches a circle. We stop computations at the time horizon $T=2.56$. In the case of the Lagrangian approach we pick $\tau = 0.002$ and the tangential redistribution parameter $\omega=1$. 
We computed evolution of the level set function in the computational domain $\Omega
\equiv \langle -4, 4 \rangle \times \langle -4, 4 \rangle$. We chose the smoothing parameter
$\epsilon=0.001$ and we did redistancing just once at $\tau_{redist}=T/2$.
Both, the initial and  final level set functions are depicted in Fig.{\ref{ellipse-test}}.

\begin{figure}
\center{
\includegraphics[width=3cm]{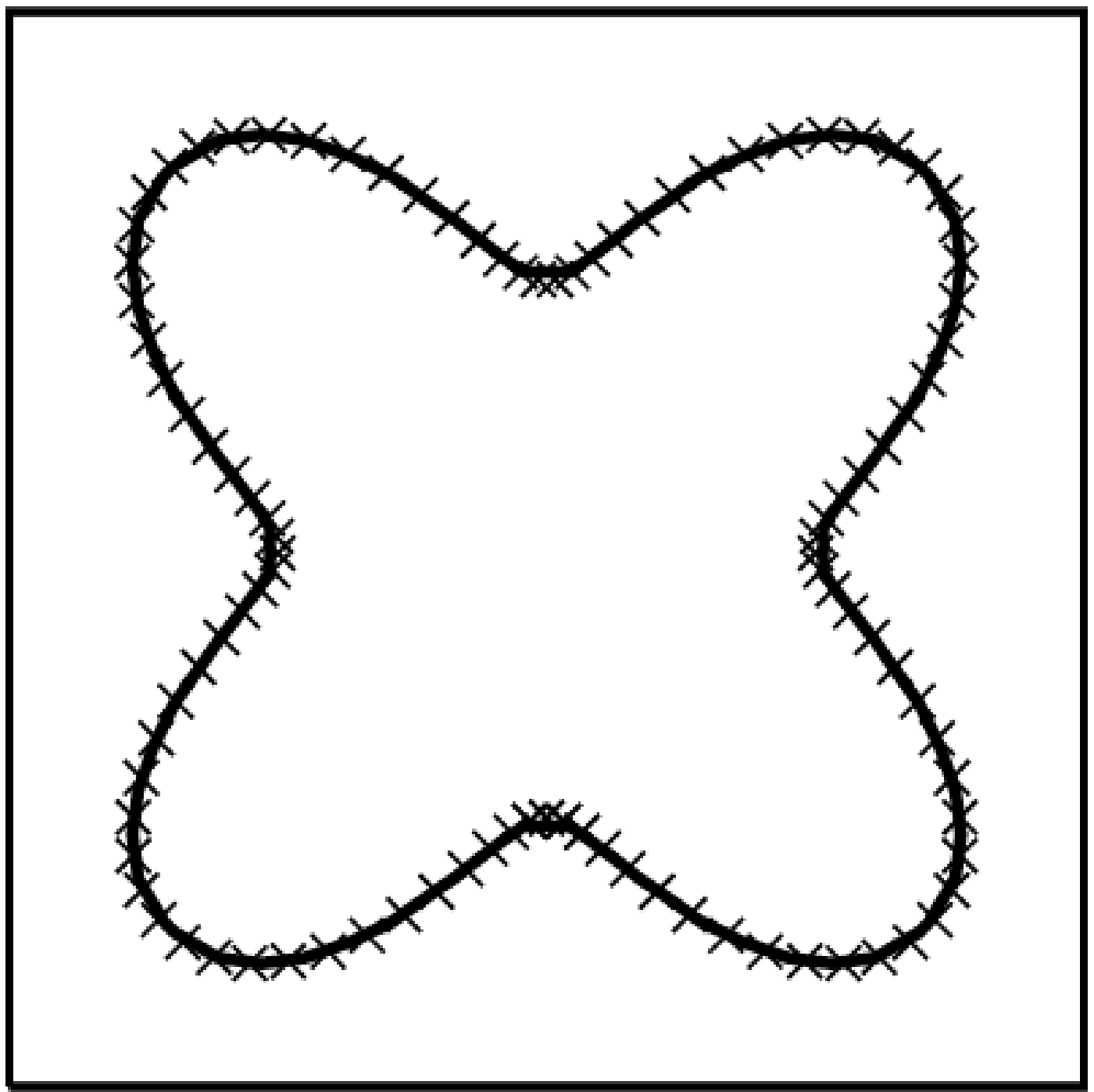}
\hglue 0.5cm
\includegraphics[width=3cm]{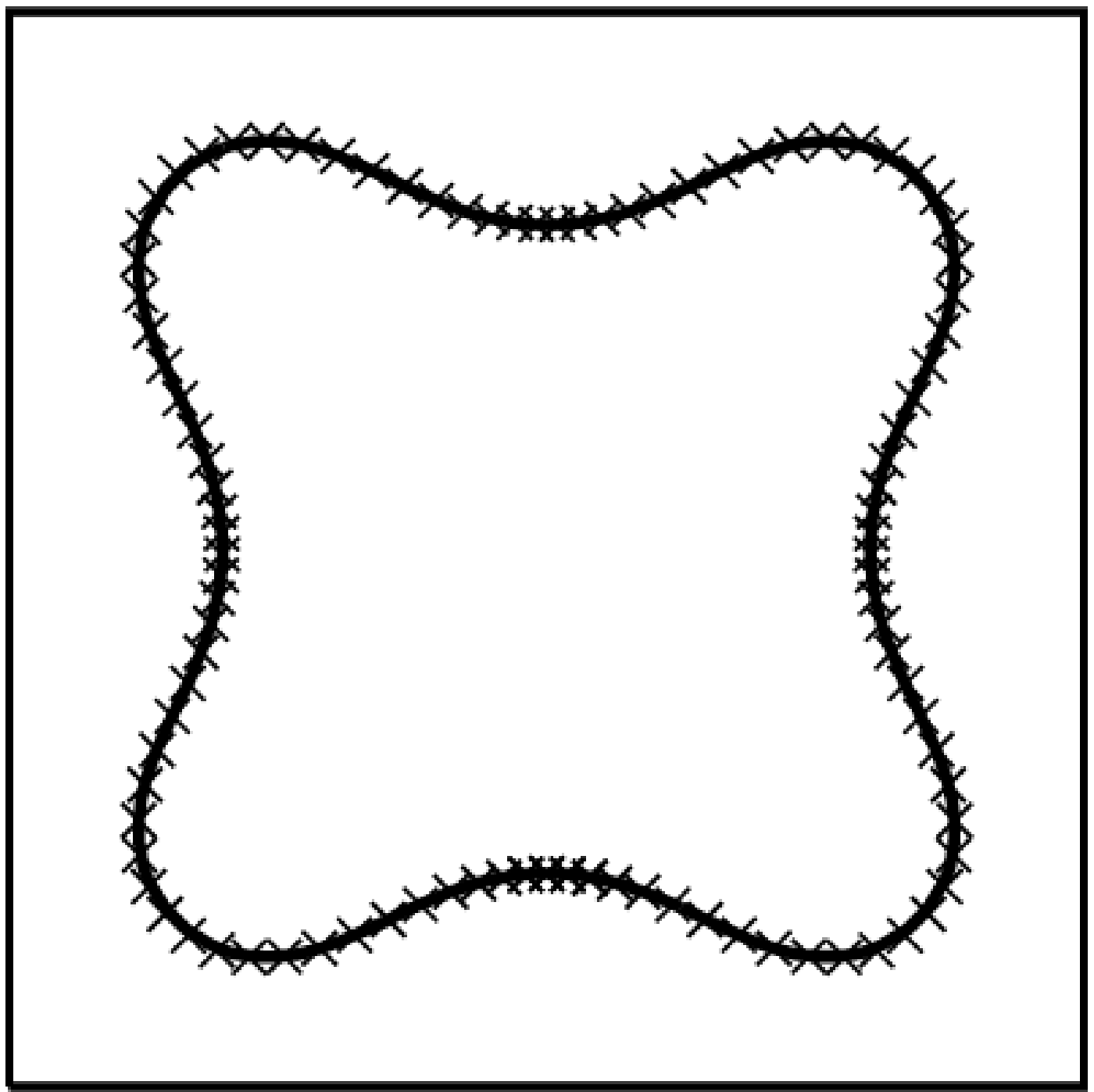}
}
\center{a) \hglue 3truecm b)}
\vskip 0.6cm
\center{
\includegraphics[width=3cm]{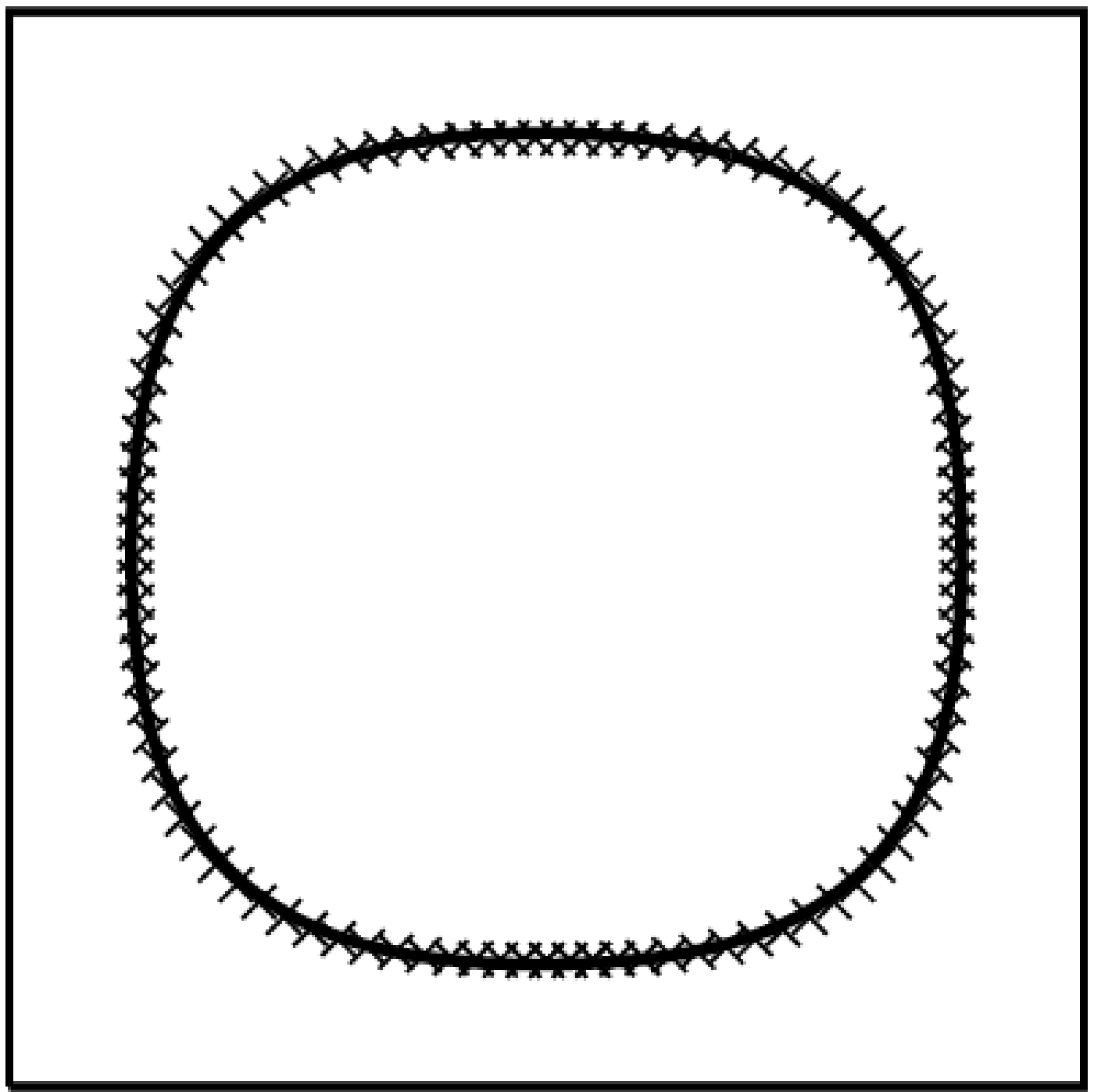}
\hglue 0.5cm
\includegraphics[width=3cm]{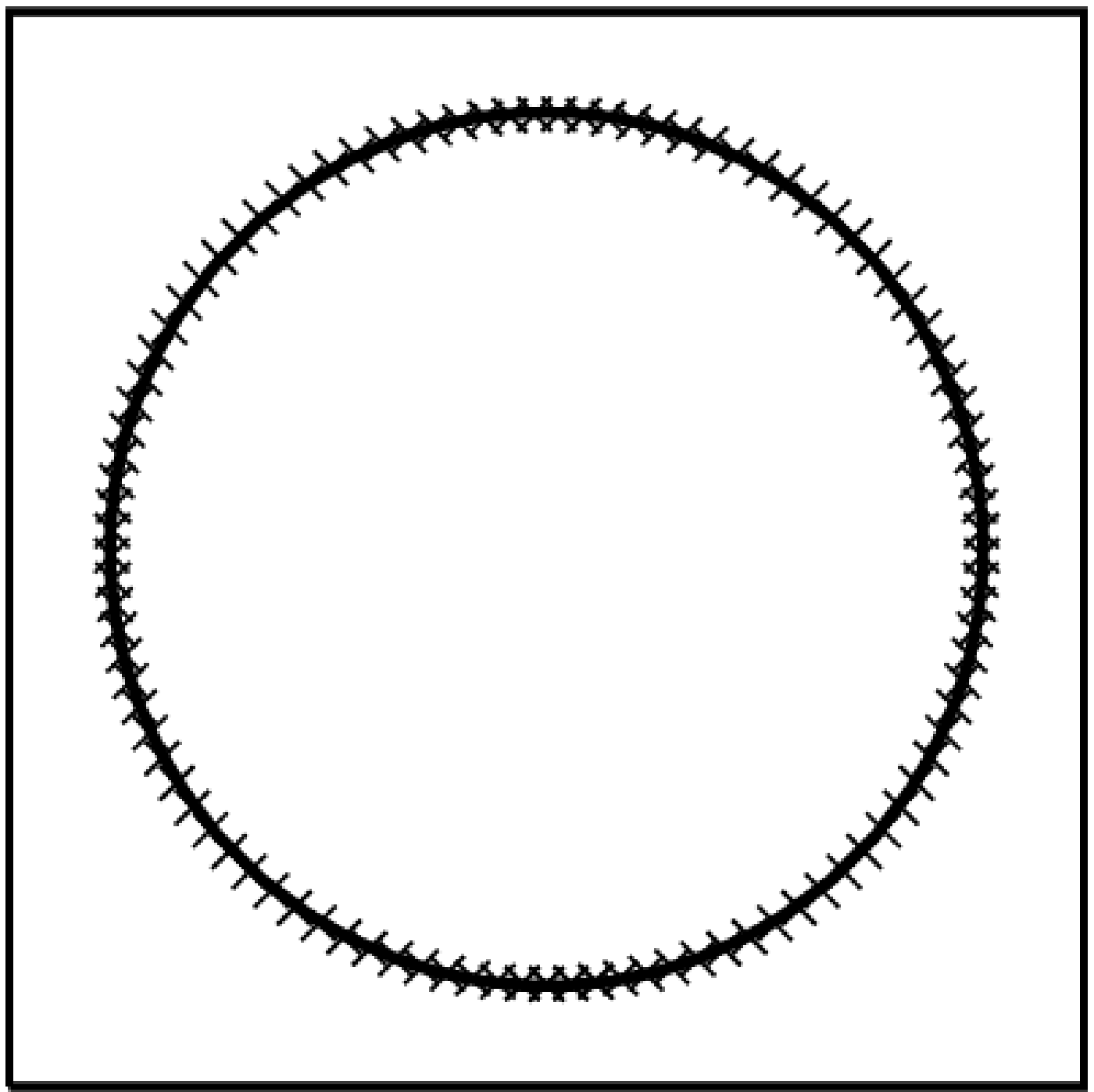}
}
\center{c) \hglue 3truecm d)}

\caption{
An initial condition a) given by (\ref{ctyrlistek-eq})  and curves computed by both methods at 
b): $t=0.001$, c): $t=0.005$ and d): $t=0.01$.}
\label{ctyrlistek-test}
\end{figure}

In Fig. \ref{ctyrlistek-test} (top) the initial condition is a non-convex  curve given by 
\begin{equation}
x^0(u)= \left( \begin{array}{c}
1 - 0.5 \cos^2( 4\pi u )\cos( 2\pi u )\\
1 - 0.5 \cos^2( 4\pi u )\sin( 2\pi u )
\end{array}
\right)
\label{ctyrlistek-eq}
\end{equation}
where $0 \leq u \leq 1$. The time evolution of such a non-convex initial curve was stopped at the time $T=0.01$. In the Lagrangian approach we picked $\tau = 10^{-5}$, and the tangential redistribution parameter was $\omega=1$. The level set function was computed in the domain 
$\Omega \equiv \langle -2, 2 \rangle
\times \langle -2, 2 \rangle$,
with $\tau = 2.5 \cdot 10^{-5}$, $\epsilon=0.001$ and $\tau_{redist}=0.001$.
Again, a comparison of the zero level set and the initial curve evolved by Lagrangian method show
compatibility of both methods in the common  time interval. In this example the 
Willmore flow quickly changes the shape of evolving curves 
from non-convex to a circular one. We show several time steps of the curve evolution
in Fig. {\ref{ctyrlistek-test}}.

\begin{figure}
\center{
\includegraphics[width=3cm]{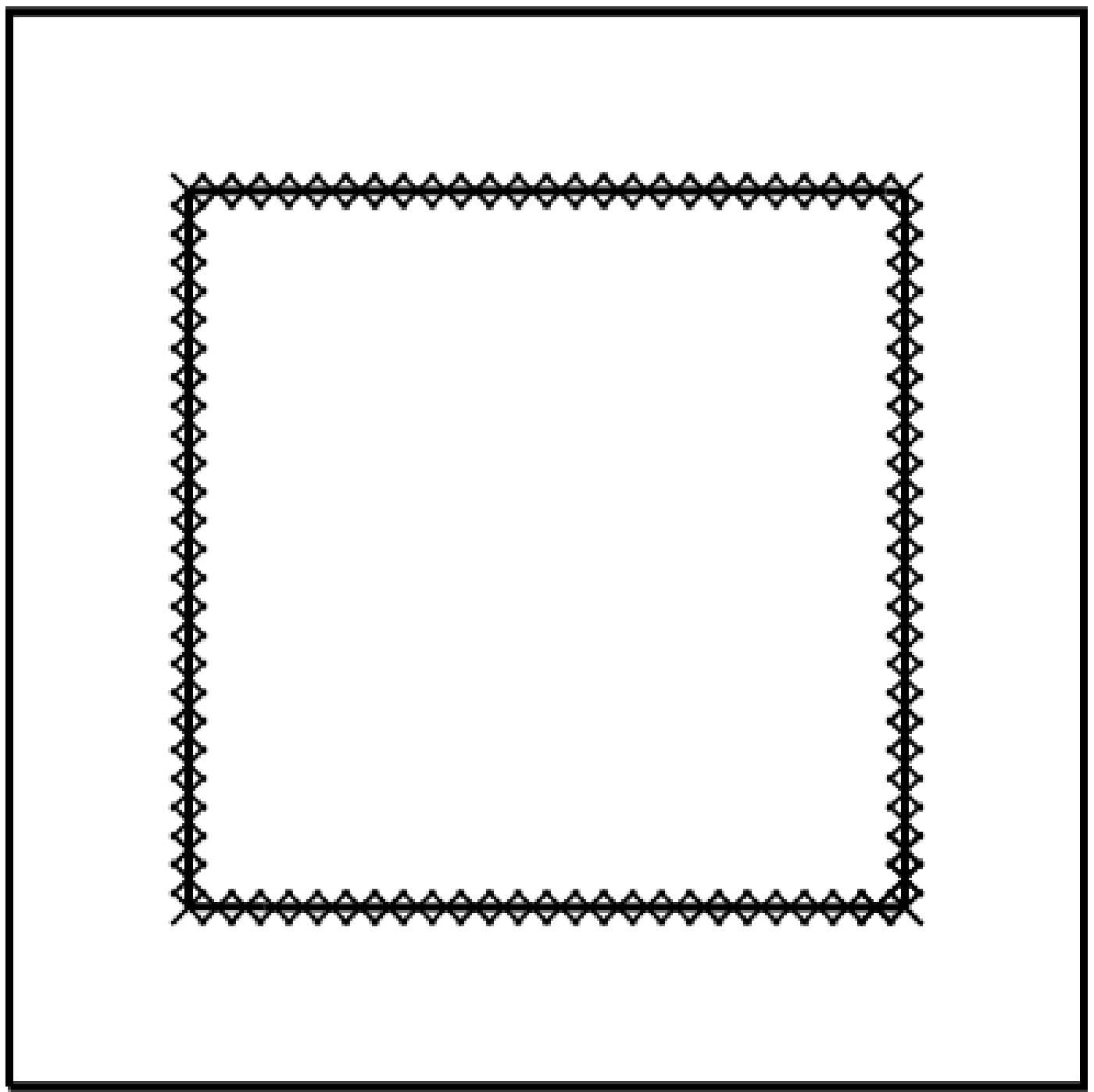}
\hglue 0.5cm
\includegraphics[width=3cm]{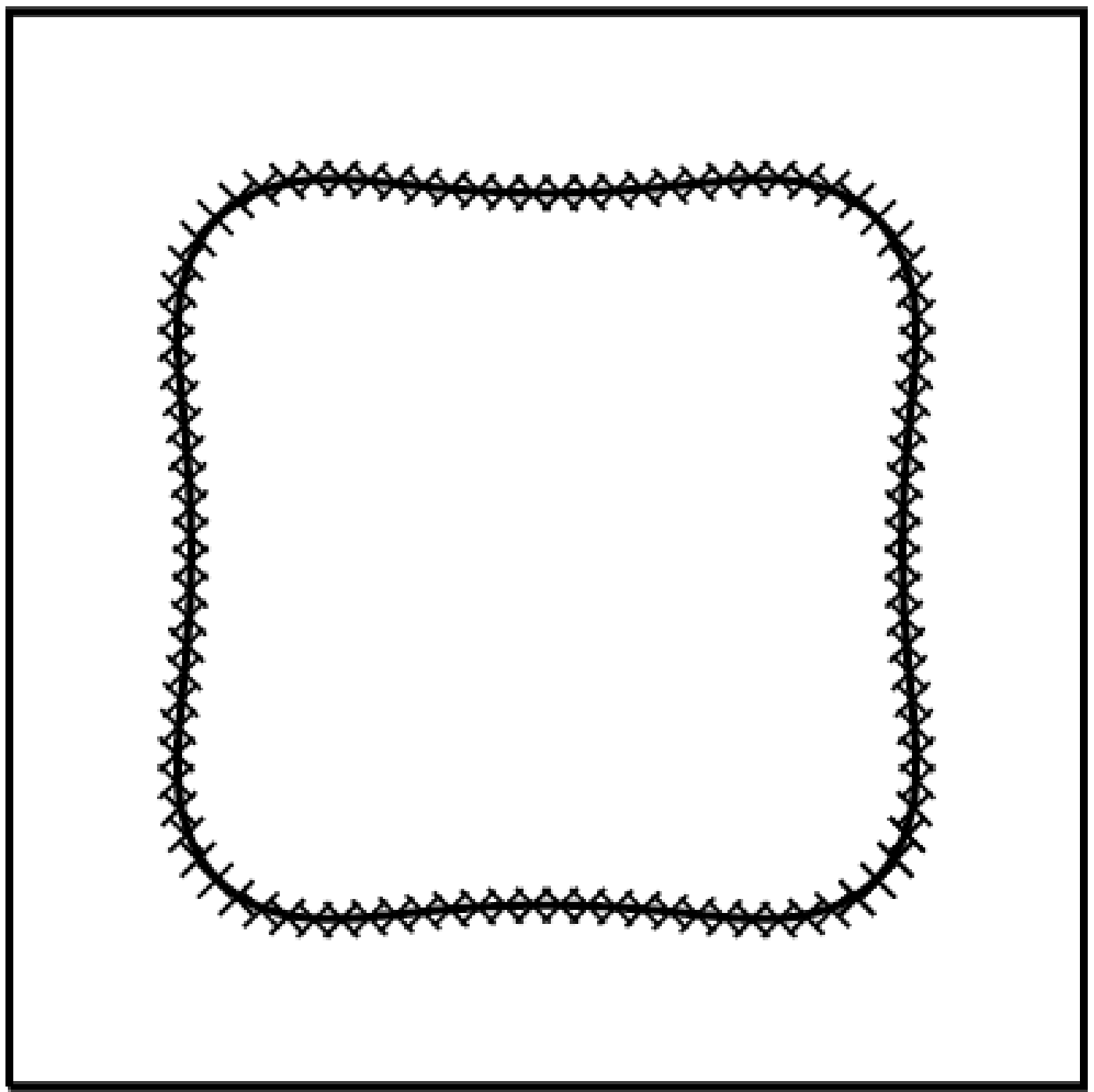}
}
\center{a) \hglue 3truecm b)}
\vskip 0.6cm
\center{
\includegraphics[width=3cm]{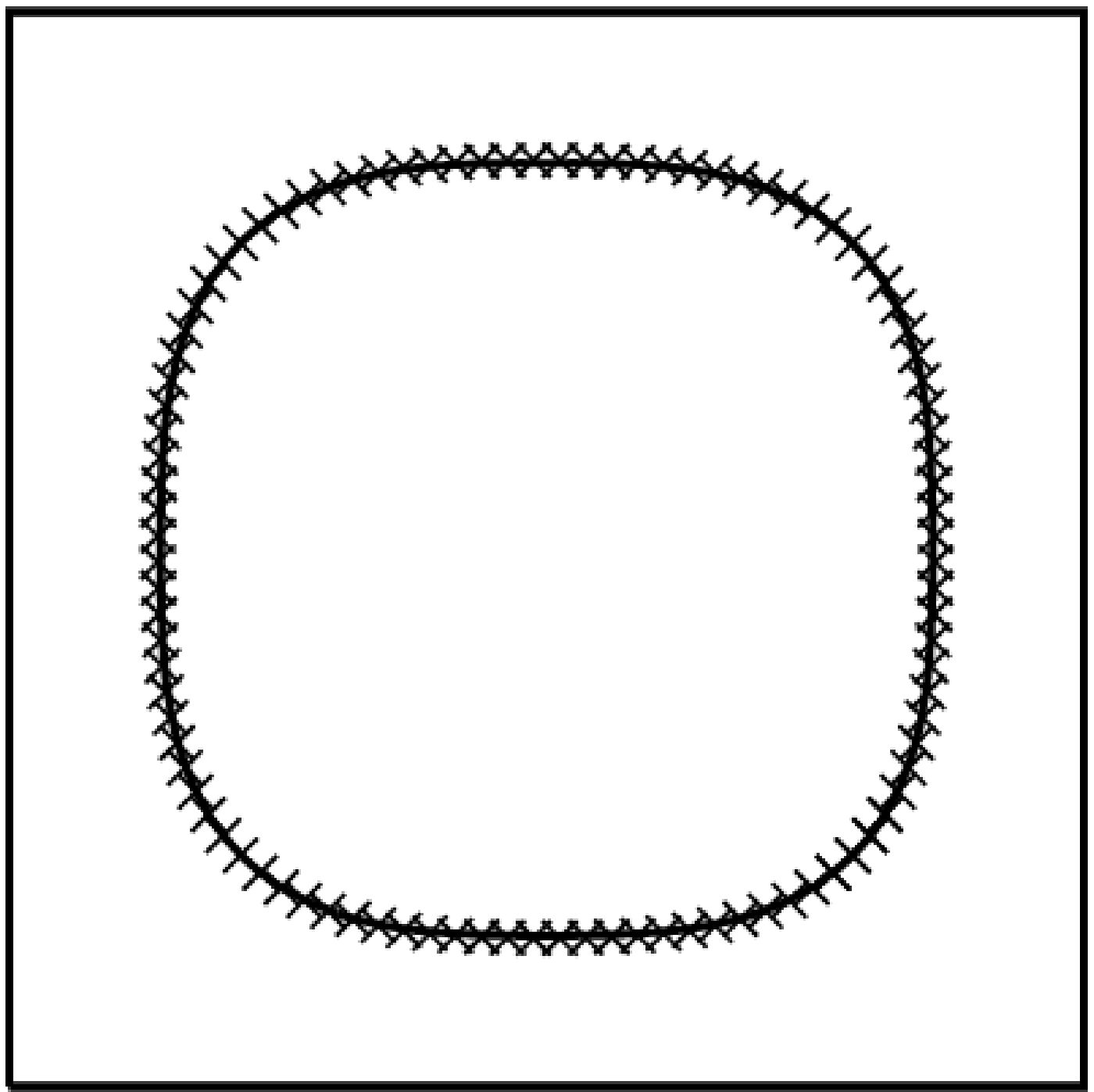}
\hglue 0.5cm
\includegraphics[width=3cm]{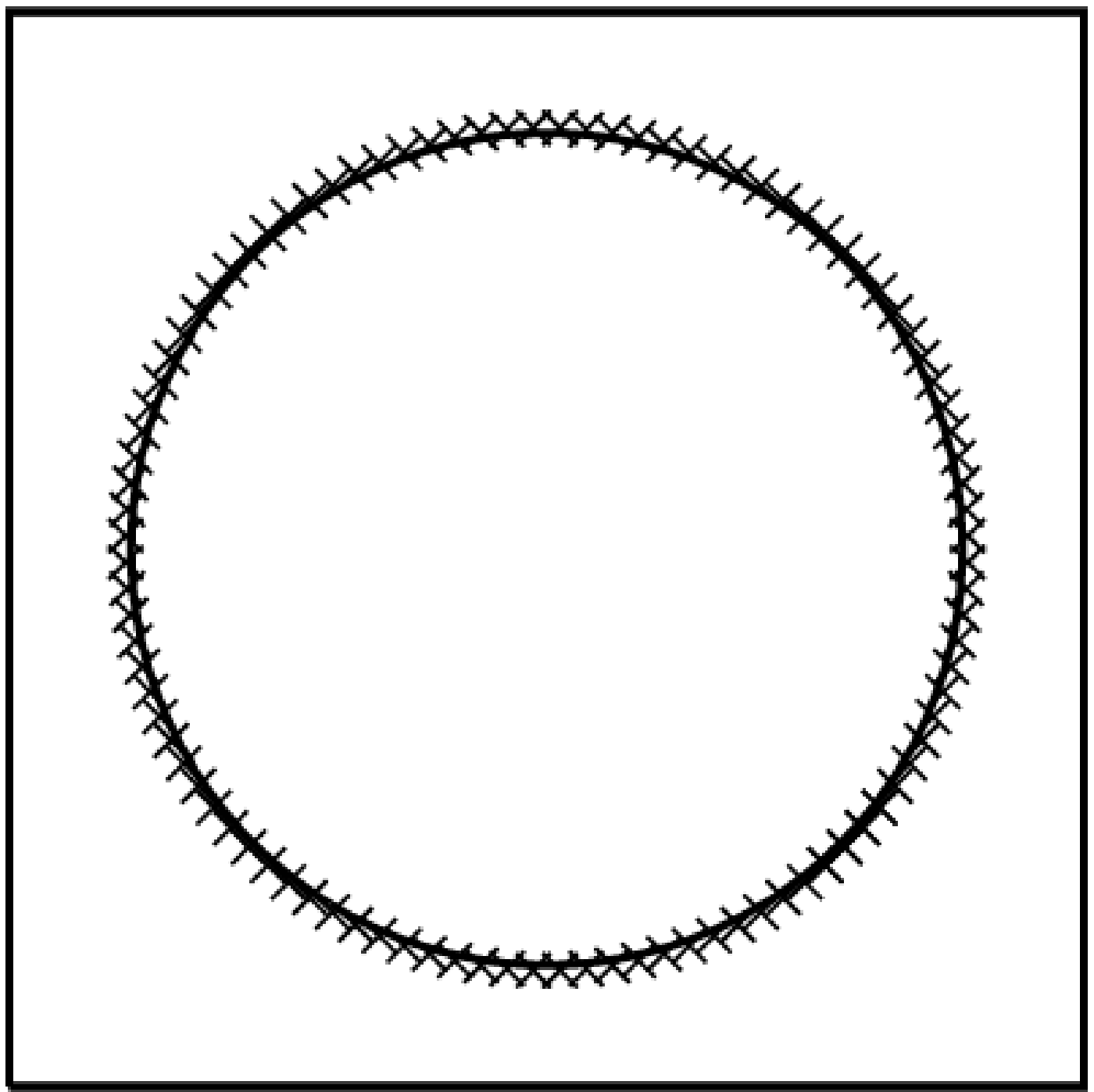}
}
\center{c) \hglue 3truecm d)}

\caption{
A square as an initial condition a) and computed evolved curves b): $t=0.001$, c): $t=0.01$
and d): $t=0.1$.}
\label{square-test}
\end{figure}

\begin{figure}
\vskip0.5cm
\center{
\includegraphics[width=4cm]{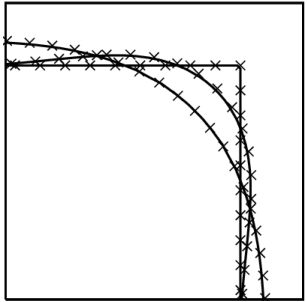}
}
\vskip 0.5cm
\caption{The detail of the square corner at the times $t=0, 0.001$ 
and $t=0.01$.}
\label{square-test-zoom}
\end{figure}

In  Fig. \ref{square-test} we present evolution with the initial curve having sharp corners (see also
a detailed close-up in  Fig. \ref{square-test-zoom}). Although the initial curve (square)  is convex, the 
evolved curve need not be convex for small times.  Concerning numerical parameters, we chose 
$n=100$ spatial nodes, $\tau=10^{-7}$ and  the tangential redistribution parameter $\omega=1$ (asymptotically uniform redistribution) in the direct Lagrangian method. As for the level-set method we took $h=0.04$, $\tau=2.5 \cdot 10^{-5}$, $\epsilon=10^{-5}$ and $\tau_{redist}=0.01$. 

\begin{figure}
\center{
\includegraphics[width=3cm]{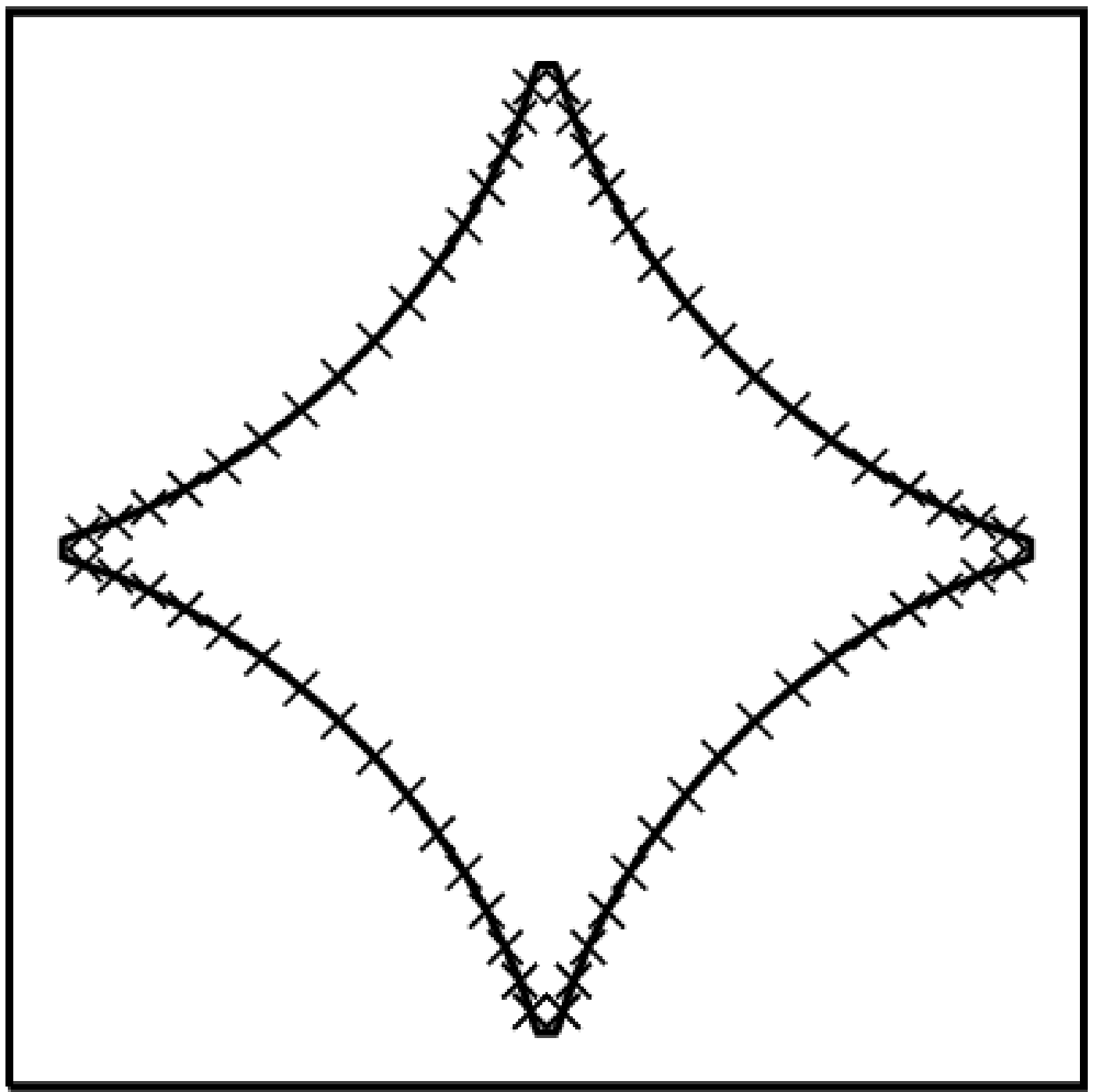}
\hglue 0.5cm
\includegraphics[width=3cm]{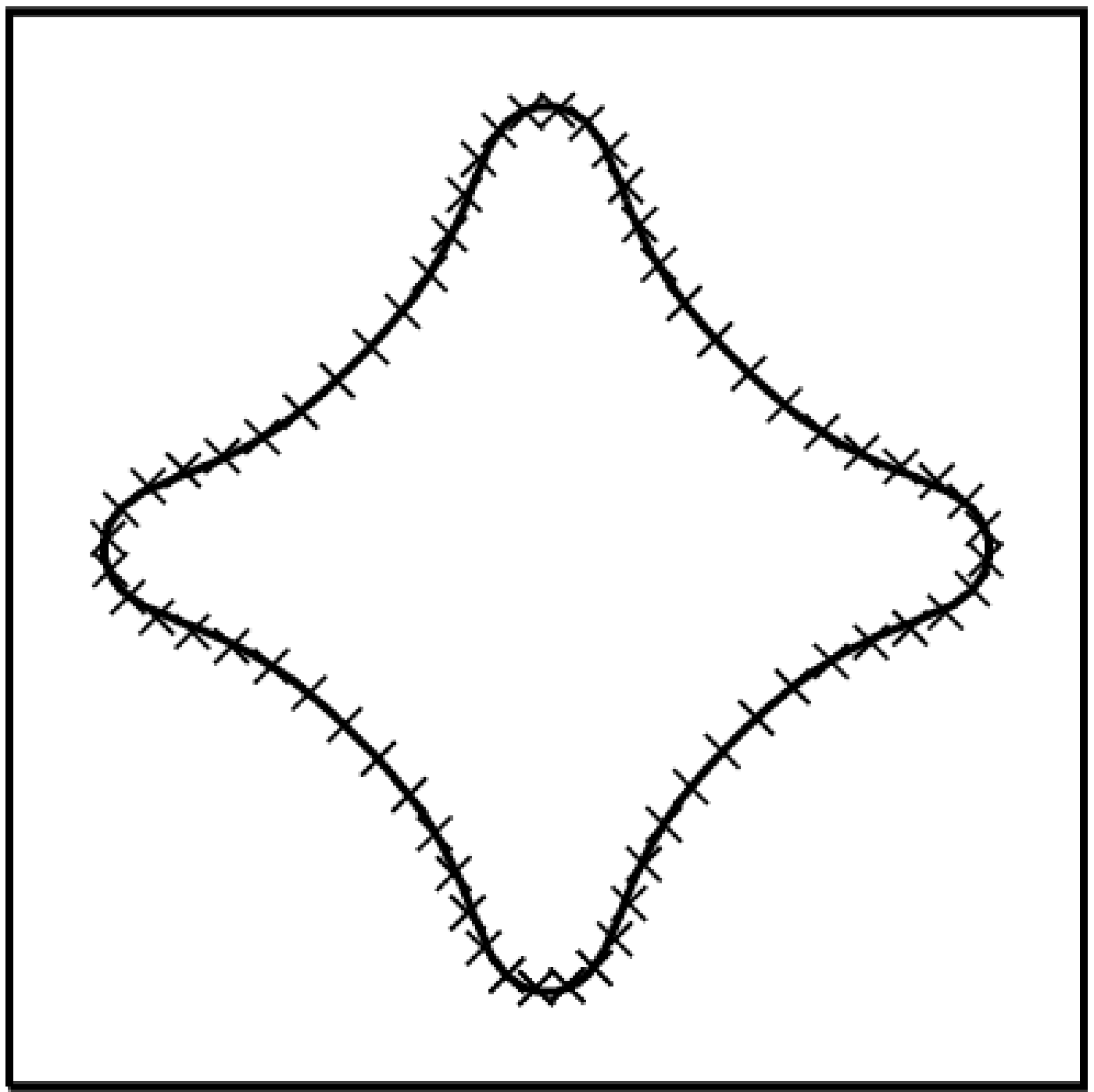}
}
\center{a) \hglue 3truecm b)}
\vskip 0.6cm
\center{
\includegraphics[width=3cm]{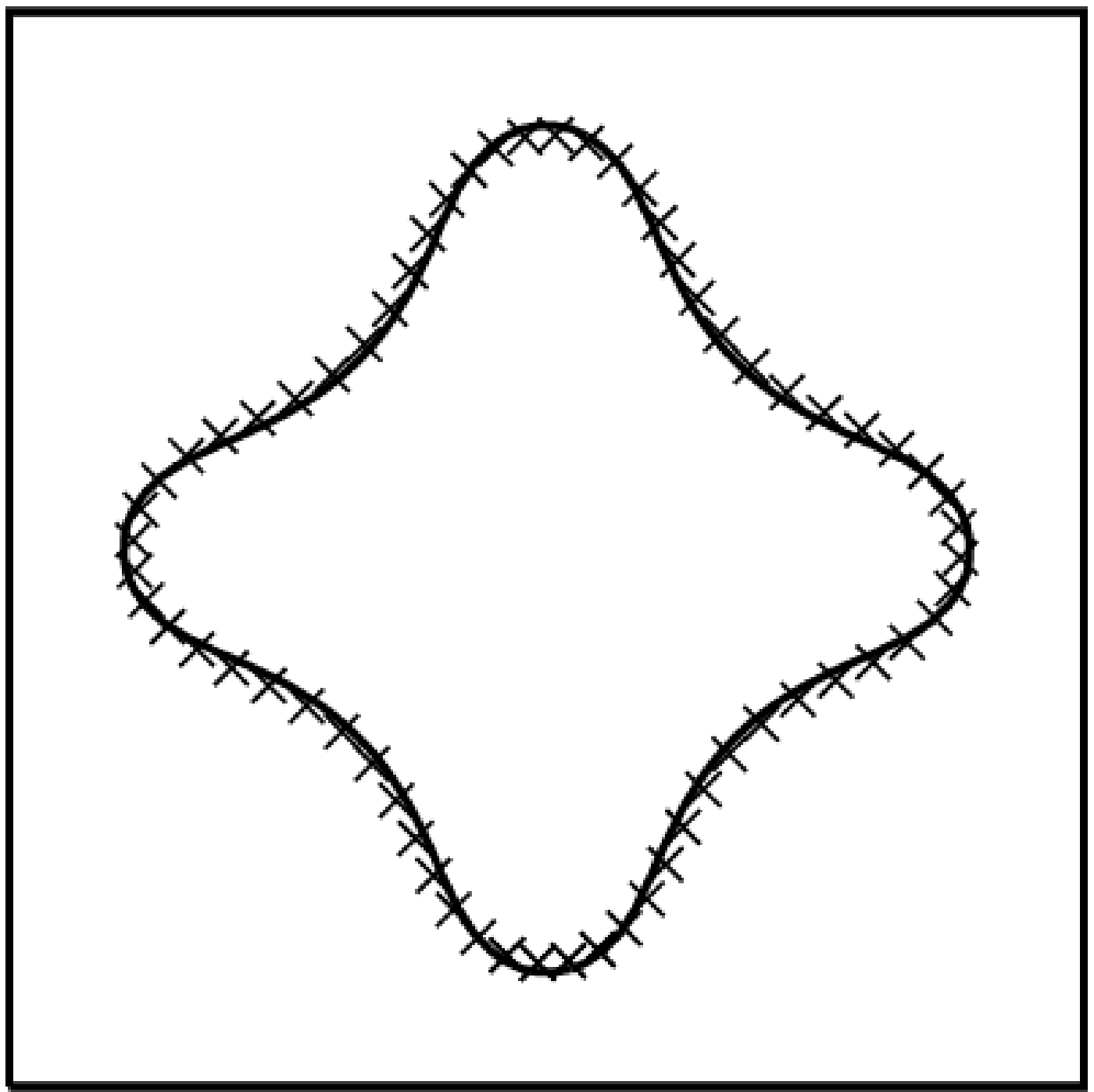}
\hglue 0.5cm
\includegraphics[width=3cm]{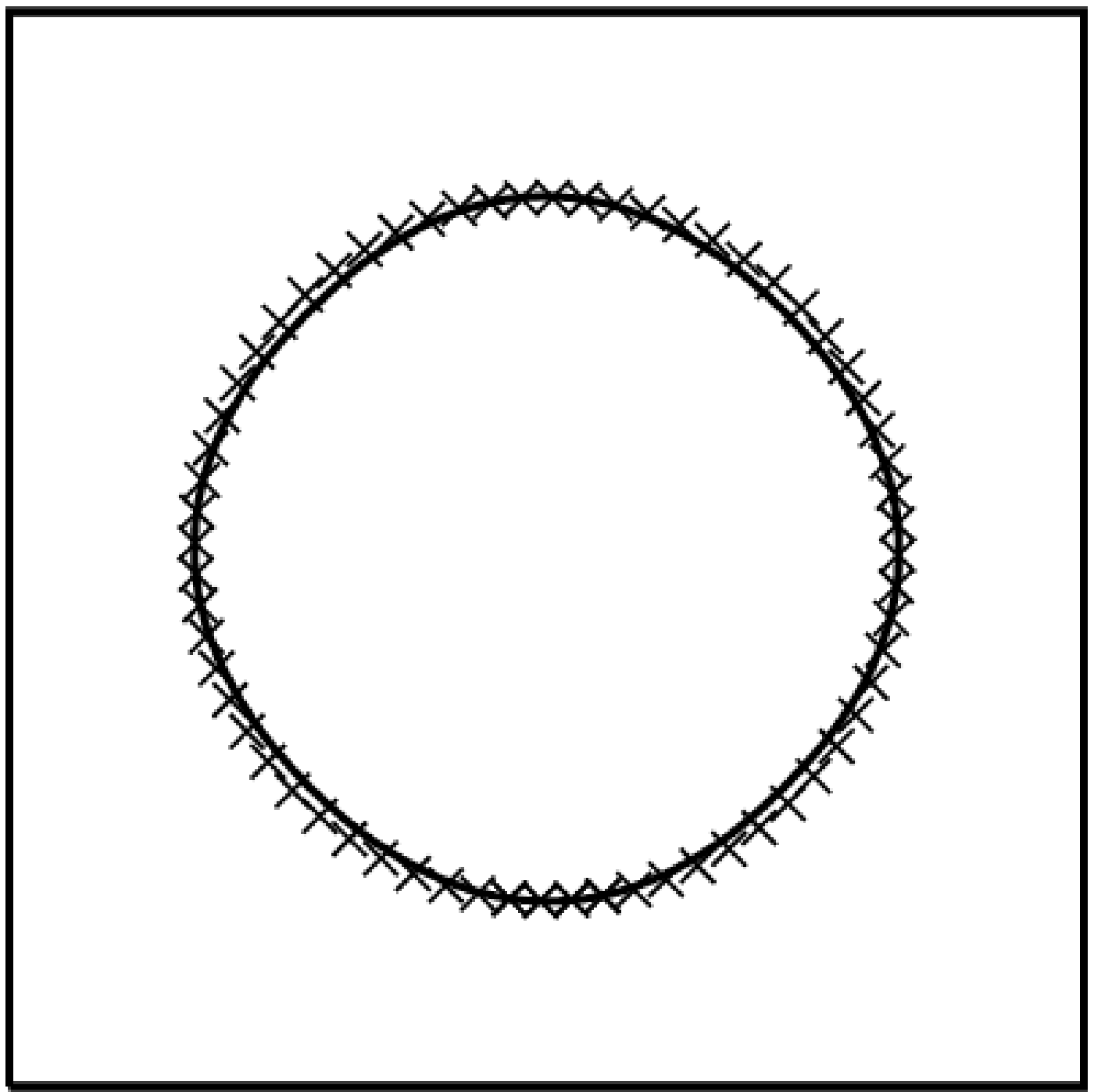}
}
\center{c) \hglue 3truecm d)}

\caption{
An asteroid as an initial condition a) and its evolution at b): $t=0.0001$; c): $t=0.0005$ and 
d): $t=0.005$.}
\label{astroid-test}
\end{figure}

Fig. \ref{astroid-test} shows comparison of the methods 
for another non-convex curve with very sharp corners. 
We chose the same numerical parameters as in the previous example 
for both the Lagrangian as well as level set methods. 
Also in this example one can observe satisfactory coincidence of numerically computed curves by both methods. 

\begin{figure}
\center{
\includegraphics[width=5.2cm]{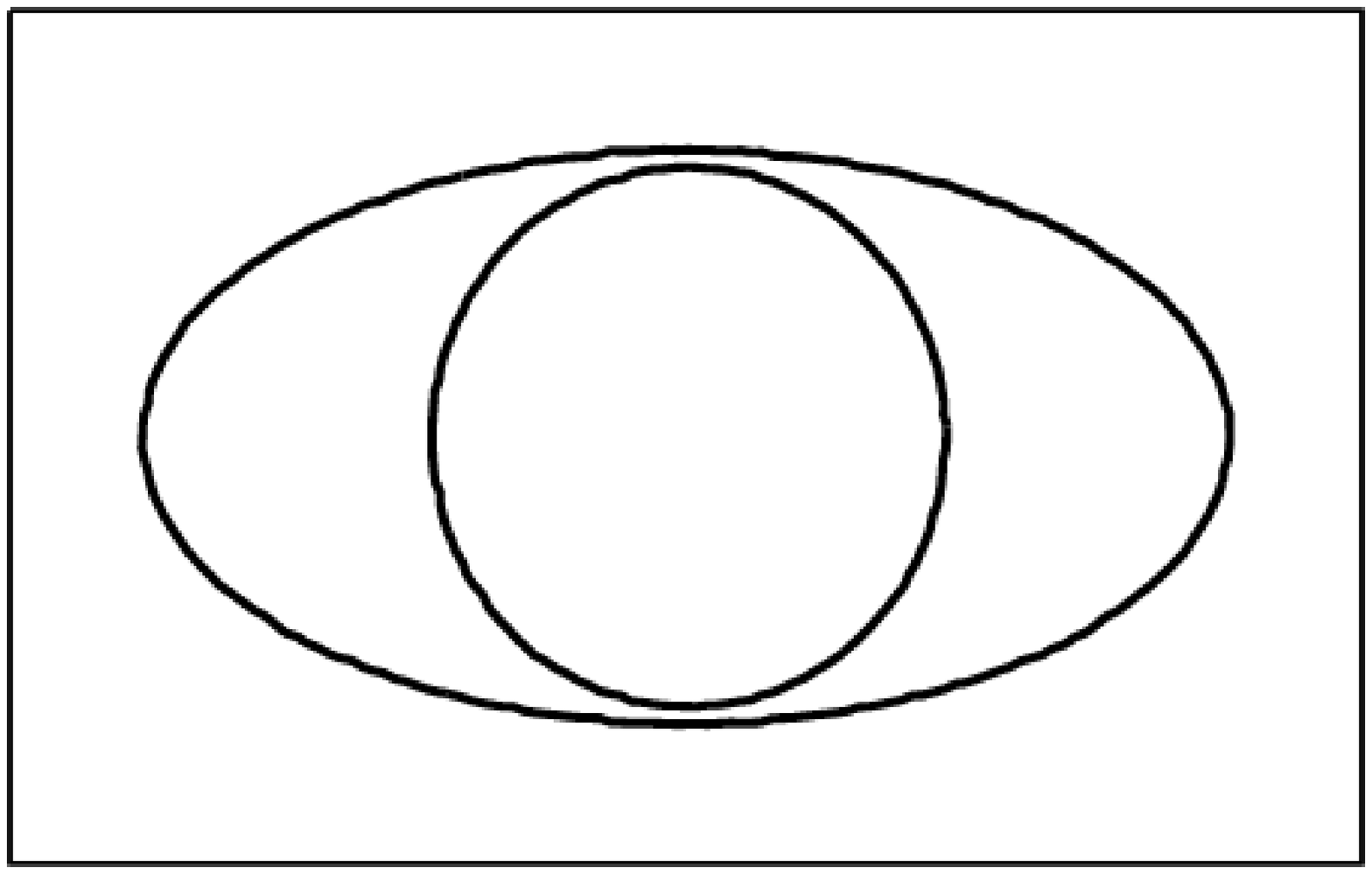}

a)

\includegraphics[width=5.2cm]{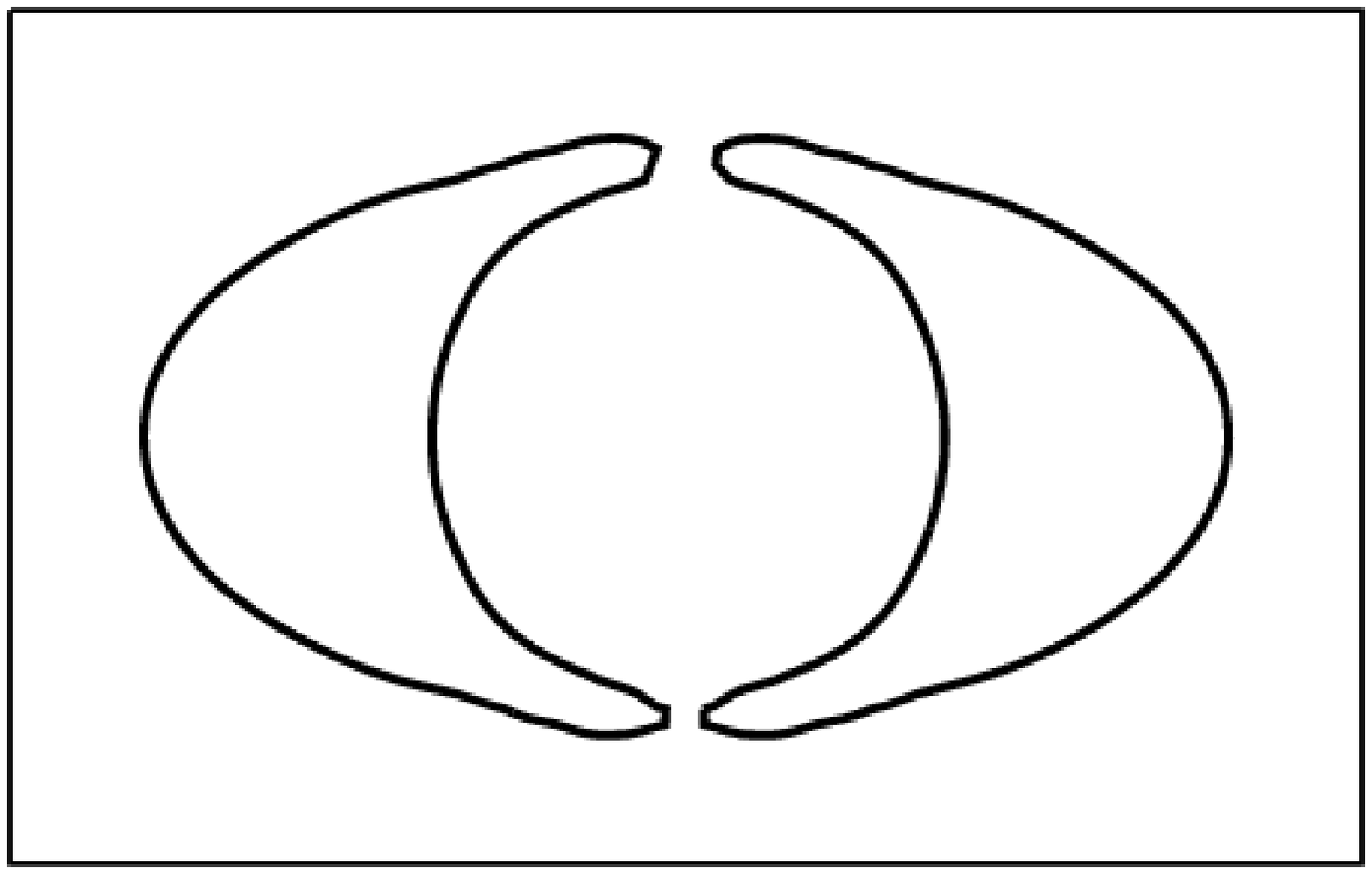}

b)

\includegraphics[width=5.2cm]{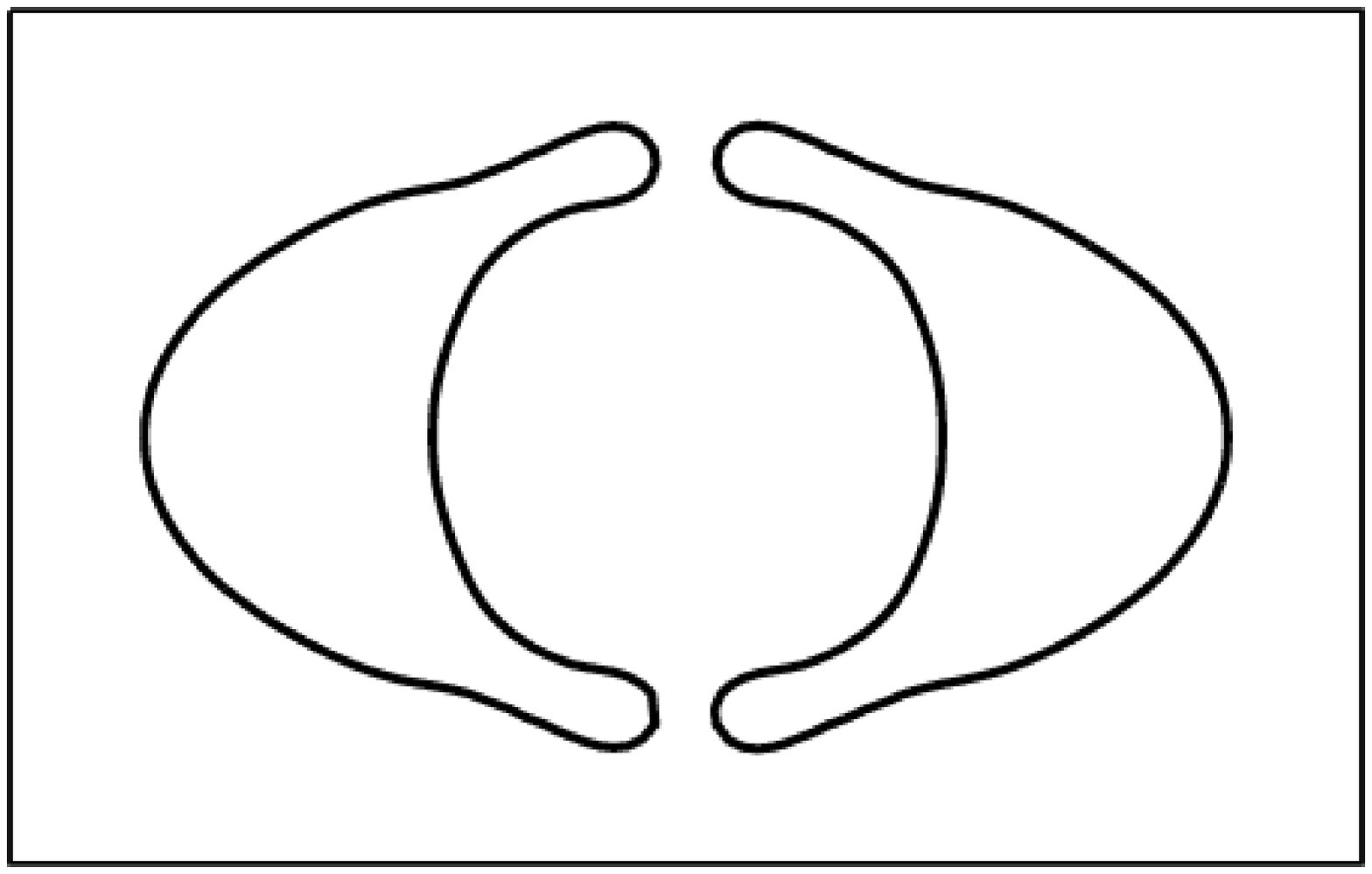}

c)

\includegraphics[width=5.2cm]{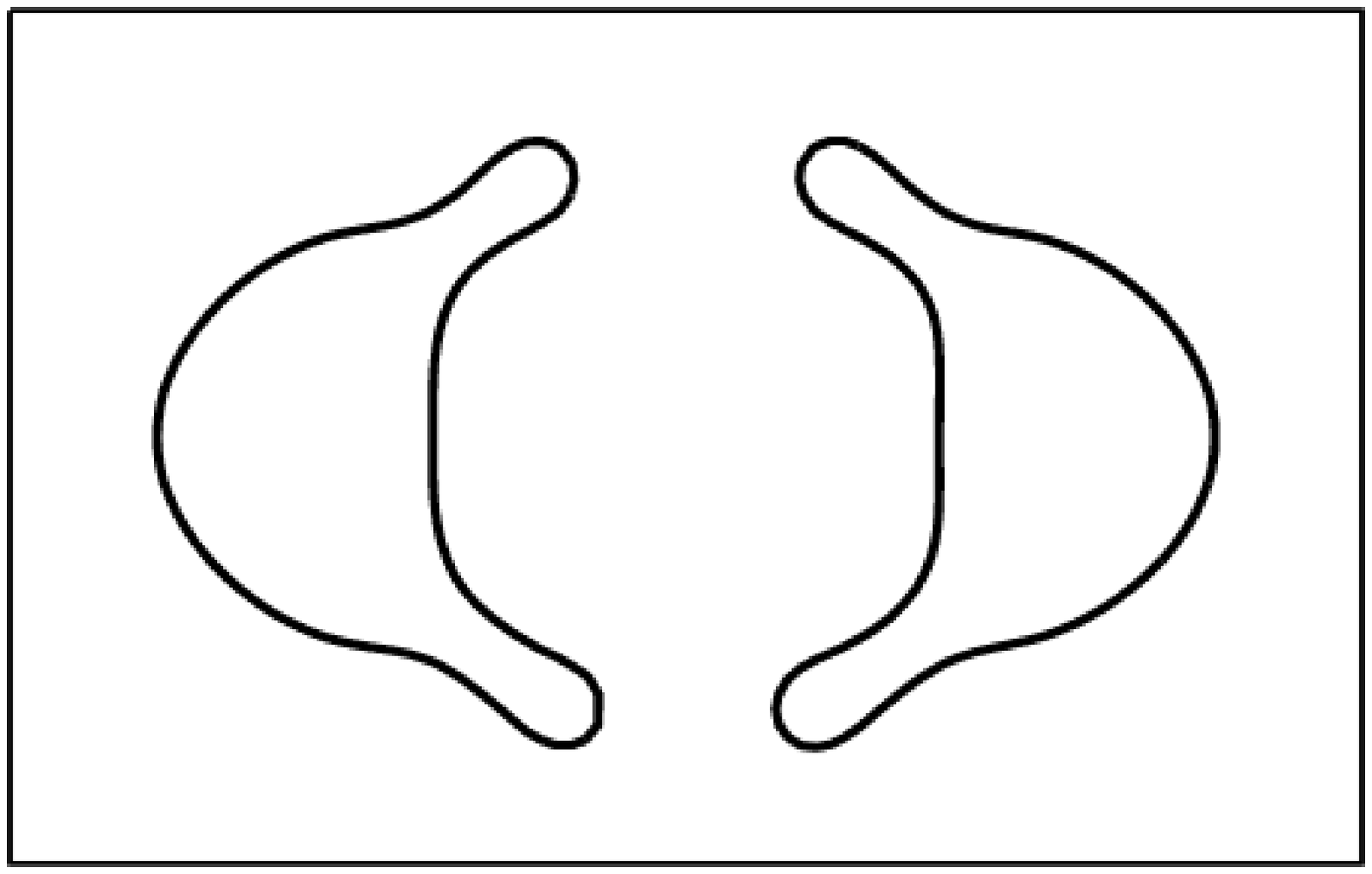}

d)

\includegraphics[width=5.2cm]{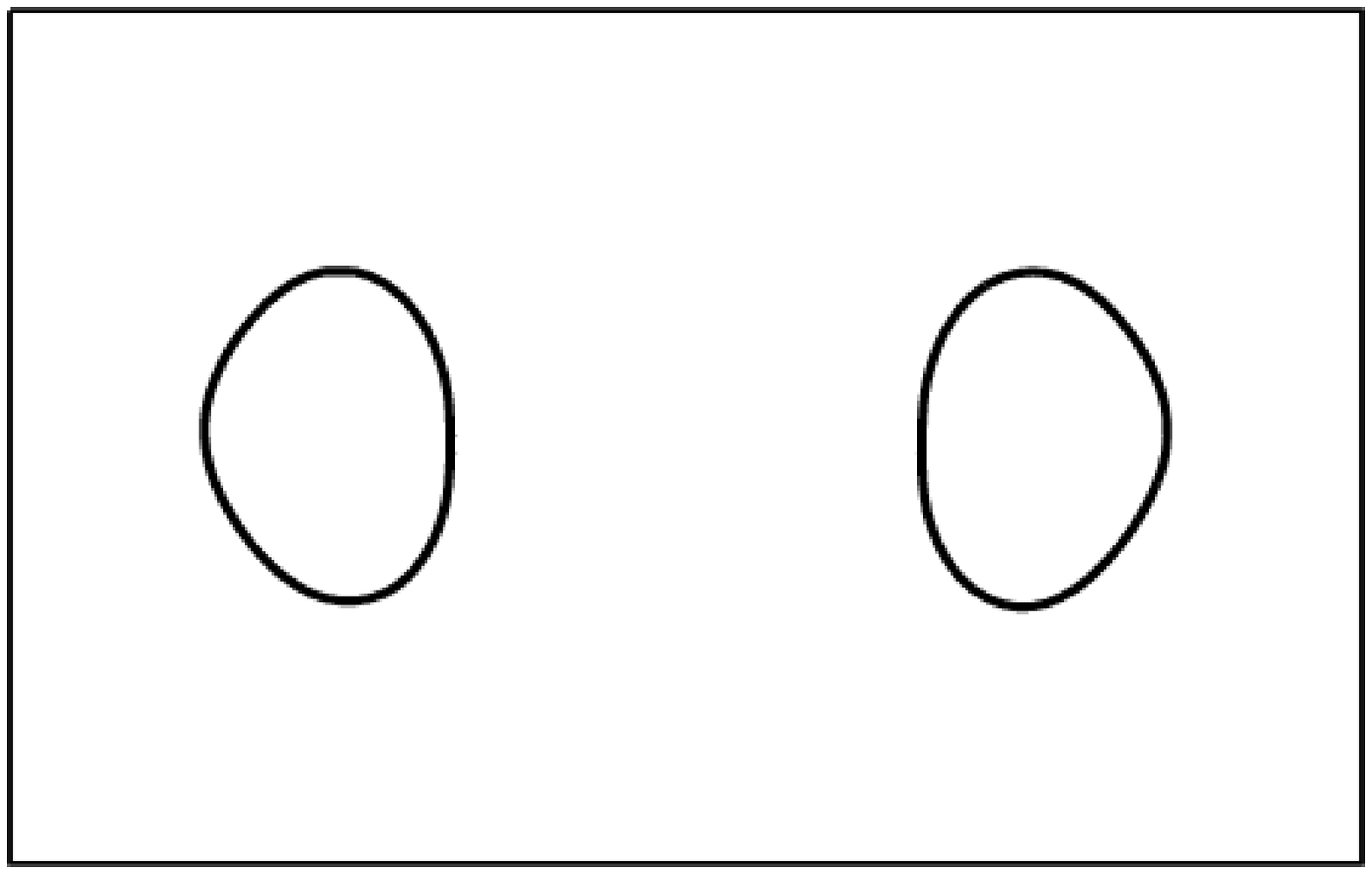}

e)
}
\caption{From top to bottom:
Initial condition composed of two curves - a circle within an ellipse a) and its evolution at b): 
$t=0.00006$; c): $t=0.0002$; d): $t=0.002$  and e): $t=0.008$.}
\label{circle-ellipse-test}
\end{figure}

Finally, in Fig. \ref{circle-ellipse-test} 
we present an example illustrating a topological change. 
It has been computed by the level set method only 
because the direct method is unable to  handle topological changes 
like pinching and splitting of curves. 
The initial zero level set consists of two almost touching curves - 
the inner curve being a circle and the outer curve being an ellipse with a shorter axis just slightly larger than the radius of the inner circle.  We then let evolve this configuration by the level set equation. The phenomenon of pinching and subsequent splitting 
of the evolved curves can be observed in this example.  
Such a behavior can be observed in the mean curvature driven 
evolution of a dumb-bell initial surface in 3D where the Grayson 
theorem does not hold.  To our best knowledge, there is no analytical proof 
of pinching-splitting phenomenon in the case of a Willmore flow of 
planar curves. Notice that 
this numerical result has been obtained only 
by using very small time steps, in the range $10^{-11} - 10^{-10}$.
Since for such small time steps we do not increase efficiency by 
using the semi-implicit scheme we use here its explicit version 
with the Runge-Kutta-Merson fourth order adaptive time solver.
As further parameters we used $h=0.0166$, $\epsilon=10^{-5}$ 
and $\tau_{redist}= 10^{-5}$.

\section{Acknowledgments}
The authors are thankful to the referee for her/his valuable comments
and suggestions that helped us to improve the final version of the paper.


\section{Appendix}

\subsection{Coefficients of the Lagrangean systems}

The coefficients for the curvature system (\ref{curvsystem}) are as follows:
\begin{eqnarray*}
a_i^j &=& \frac 1{q_{i-1}^j r_{i-1}^j q_{i-2}^j},\ \ 
e_i^j = \frac 1{q_i^j r_{i+1}^j q_{i+1}^j},\\
b_i^j &=& -\left(\frac 1{r_i^j q_i^j q_{i-1}^j}
+\frac 1{r_i^j (q_{i-1}^j)^2}+\frac 1{(q_{i-1}^j)^2 r_{i-1}^j}\right.\\
&&\left. +\ \ \ \ \frac 1{q_{i-1}^j r_{i-1}^j
q_{i-2}^j}\right)+\frac{\alpha_{i-1}^j}{2}\,,
\\
d_i^j &=& -\left(\frac 1{q_i^j r_{i+1}^j q_{i+1}^j} + 
\frac 1{(q_i^j)^2 r_{i+1}^j} +
\frac 1{r_i^j (q_i^j)^2}\right.\\
&&\left.+\ \ \ \ \frac 1{r_i^j q_i^j q_{i-1}^j}\right)
-\frac{\alpha_{i}^j}{2}\,,
\\
c_i^j &=& \frac 1{(q_i^j)^2 r_{i+1}^j} +
\frac 1{r_i^j (q_i^j)^2}+\frac 2{r_i^j q_i^j q_{i-1}^j}
+\frac 1{r_i^j (q_{i-1}^j)^2}\\
&+&\frac 1{(q_{i-1}^j)^2 r_{i-1}^j}+\frac{r_i^j}{\tau}
- r_i^{j-1} k_i^{j-1}\beta_i^{j-1}+
\frac{\alpha_i^j}{2}-\frac{\alpha_{i-1}^j}{2}\,,
\\
\\
f_i^j &=& \frac {r_i^j}{\tau}
k_i^{j-1}+\frac{(k_{i}^{j-1})^3 -(k_{i-1}^{j-1})^3}
   {2 q_{i-1}^j}-\frac{(k_{i+1}^{j-1})^3-(k_{i}^{j-1})^3}
   {2 q_i^j}
\end{eqnarray*}
where we used
following approximation of third order derivative terms on
boundaries of flowing finite volume in (\ref{curvdiscr}):
\begin{eqnarray*}
&&\partial_s^3 k(x_i)-\partial_s^3
k(x_{i-1})\approx \\
\\
&&\approx 
\frac{\partial_s^2 k(\tilde x_{i+1})-\partial_s^2 k(\tilde x_{i})}{q_i}
-\frac{\partial_s^2 k(\tilde x_{i})-\partial_s^2 k(\tilde x_{i-1})}{q_{i-1}}
\\
&&\approx\dots\approx 
\frac 1{q_i r_{i+1} q_{i+1}} k_{i+2} 
+\frac 1{q_{i-1} r_{i-1} q_{i-2}} k_{i-2}-\\
&&\left(\frac 1{q_i r_{i+1} q_{i+1}} + \frac 1{q_i^2 r_{i+1}} +
\frac 1{r_i q_i^2}+\frac 1{r_i q_i q_{i-1}}\right) k_{i+1} +
\\
&&\left(\frac 1{q_i^2 r_{i+1}} +      
\frac 1{r_i q_i^2}+\frac 2{r_i q_i q_{i-1}}
+\frac 1{r_i q_{i-1}^2}+\frac 1{q_{i-1}^2 r_{i-1}}\right) k_i -
\\
&&\left(\frac 1{r_i q_i q_{i-1}}
+\frac 1{r_i q_{i-1}^2}+\frac 1{q_{i-1}^2 r_{i-1}}
+\frac 1{q_{i-1} r_{i-1} q_{i-2}}\right) k_{i-1}. \\
\end{eqnarray*}
Using a similar strategy for approximation of the third order derivatives
of position vector on boundaries of flowing dual volume we can write coefficients
of (\ref{systempos}):
\begin{eqnarray}
{\cal A}_i^j &=& \frac 1{r_{i}^j q_{i-1}^j r_{i-1}^j},\ \ \  
{\cal C}_i^j = \frac{q_i^j}{\tau}
-({\cal A}_i^j+{\cal B}_i^j+{\cal D}_i^j+{\cal E}_i^j),\nonumber \\ 
{\cal E}_i^j &=& \frac 1{r_{i+1}^j q_{i+1}^j r_{i+2}^j},\ \ \ \  
{\cal F}_i^j = \frac {q_i^j}{\tau}x_i^{j-1},\nonumber\\
{\cal B}_i^j &=& -\left(\frac 1{r_i^j q_{i-1}^j r_{i-1}^j}
+\frac 1{(r_i^j)^2 q_{i-1}^j}+\frac 1{(r_i^j)^2 q_i^j}\right.
\nonumber\\
&&\left.+\ \ \ \ \frac 1{r_i^j q_i^j r_{i+1}^j}\right) 
+\frac 32 \frac{(k_i^j)^2}{r_i^j} + \frac{\alpha_{i}^j}{2}\,,
\nonumber\\
{\cal D}_i^j &=& -\left(\frac 1{r_i^j q_{i}^j r_{i+1}^j}
+\frac 1{(r_{i+1}^j)^2 q_{i}^j}+\frac 1{(r_{i+1}^j)^2 q_{i+1}^j}\right.
\nonumber\\
&&\left.+\ \ \ \ \frac 1{r_{i+1}^j q_{i+1}^j r_{i+2}^j}\right)
+\frac 32 \frac{(k_{i+1}^j)^2}{r_{i+1}^j} 
-\frac{\alpha_{i}^j}{2}\,.
\nonumber
\end{eqnarray}

\subsection{Coefficients of the level set system}

The coefficients $A_{ij}^{rs}$ of the 21-diagonal system matrix
(\ref{systemlevelset}) are given by:

\begin{eqnarray*}
A_{ij}^{00} &=& 1 + \frac{\tau Q_{ij}}{h^4} \sum_{r,s \in \left\{ -1, 1\right\}} \left[
   \frac{h^2\left( \hat w_{ij}^{rs;n-1}\right)^2}{2 \left( Q_{ij}^{rs;n-1}\right)^3} \right.+ \\
   &&{E}_{1+|s|,2-|r|;ij}^{r,s}
      \left( \frac{\bar Q_{i+r,j+s}^{n-1}}{Q_{i+r,j+s}^{-r,-s;n-1}} + \frac{\bar Q_{ij}^{n-1}}{Q^{\ast;n-1}_{ij}} \right) +  \\
   &&\left. \frac14 {E}_{1+|s|,1+|r|;ij}^{r,s} 
      \left( \frac{\bar Q_{i+s,j+r}^{n-1}}{Q_{i+s,j+r}^{-r,-s;n-1}} - \frac{\bar Q_{i-s,j-r}^{n-1}}{Q_{i-s,j-r}^{sr;n-1}} \right)
\right],
\end{eqnarray*}

\begin{eqnarray*}
A_{ij}^{rs} &=& \frac{\tau}{h^4} \bar Q^{n-1}_{ij} \left[ 
- {E}_{1+|s|, 2 - |r|;ij}^{rs; n-1}
   \left(
   \frac{\bar Q^{n-1}_{i+r,j+s} }{Q^{\ast;n-1}_{i+r,j+s}} + 
   \frac{\bar Q^{n-1}_{ij}}{Q_{ij}^{rs; n-1}} 
   \right) \right.\\
   &+& 
\frac14 {E}_{1+|s|,1+|r|;ij}^{rs;n-1}
   \left(
   \frac{\bar Q_{i+r+s,j+r+s}^{n-1}}{Q_{i+r+s,j+r+s}^{-s,-r; n-1}} -
   \frac{\bar Q_{i+r-s,i-r+s}^{n-1}}{Q_{i+r-s,j-r+s}^{sr; n-1} }
   \right) \\
   &+& 
\frac14 {E}_{1+|r|,1+|s|;ij}^{sr; n-1}
   \left(
   \frac{\bar Q_{i+r+s,j+r+s}^{n-1}}{Q_{i+r+s,j+r+s}^{-s,-r; n-1}} -
   \frac{\bar Q_{i+r,j+s}^{n-1}}{Q^{\ast; n-1}_{i+r,j+s}}
   \right) \\
   &-& 
\frac 14 {E}_{1+|r|,1+|s|;ij}^{-s,-r; n-1} 
   \left(
   \frac{\bar Q_{i+r-s,j-r+s}^{n-1}}{Q_{i+r-s,j-r+s}^{sr; n-1}}-
   \frac{\bar Q_{i+r,j+s}^{n-1}}{Q^{\ast; n-1}_{i+r,j+s}}
   \right) \\
   &-& \left( 
   {E}_{1+|s|,2-|r|;ij}^{-r,-s; n-1} +
   {E}_{1+|r|,2-|s|;ij}^{sr; n-1} \right.\\
&& +   \left. {E}_{1+|r|,2-|s|}^{-s,-r; n-1} \right) \frac{\bar
Q_{ij}^{n-1}}{Q_{ij}^{rs; n-1}}
-\left. \frac{ h^2 \left( \hat w_{ij}^{rs;n-1}\right)^2}{2 (Q_{ij}^{rs; n-1})^3 }
\right],
\end{eqnarray*}
for $|r|+|s|=1$. Here we have denoted by $Q^{\ast;n}_{i,j}$ the harmonic
average of $Q^{rs; n}_{ij}$ defined as:
\[
\frac{1}{Q^{\ast;n}_{i,j}} = \sum_{|r|+ |s| = 1} \frac{1}{Q^{rs; n}_{ij}}\,.
\]
For $|r| = 1$ and $|s| = 1$ we have the expression:
\begin{eqnarray*}
A_{ij}^{rs} &=& rs \frac{\tau}{h^4} \bar Q_{ij}^{n-1} \left[
\frac{{E}_{11;ij}^{r0; n-1} \bar Q_{i+r,j}^{n-1}}{Q_{i+r,j}^{0s;n-1}} +
\frac{{E}_{22;ij}^{0s;n-1} \bar Q_{i,j+s}^{n-1}}{Q_{i,j+s}^{r0; n-1}}
\nonumber \right. \\ 
&+&  \frac14 {E}_{12;ij}^{r0; n-1} 
\left( \frac{\bar Q_{i,j+s}^{n-1}}{Q_{i,j+s}^{r0; n-1}} - 
\frac{\bar Q_{i+r,j+s}^{n-1}}{Q^{\ast; n-1}_{i+r,j+s}}
\right) \\
&+&  
\frac14 {E}_{21;ij}^{0s; n-1}
\left(
\frac{\bar Q_{i+r,j}^{n-1}}{Q_{i+r,j}^{0s; n-1}} - 
\frac{Q_{i+r,j+s}}{Q^{\ast; n-1}_{i+r,j+s}}
\right) 
\\
&-& 
\left.  \frac14 \frac{{E}_{12;ij}^{-r,0; n-1} \bar
Q_{i,j+s}^{n-1}}{Q_{i,j+s}^{r,0; n-1}} 
- \frac14 \frac{{E}_{21;ij}^{0,-s; n-1} \bar Q_{i+r,j}^{n-1}}{Q_{i+r,j}^{0,s;
n-1}} \right]\,.
\end{eqnarray*}
Next, for $|r|=2$ or $|s|=2$ such that $|r| + |s| = 2$ we have
\begin{eqnarray*}
A_{ij}^{rs} &=& \frac{\tau}{h^4} 
\bar Q_{ij}^{n-1} \left(
{E}_{1+|\bar s|,1+|\bar s|;ij}^{\bar r,\bar s; n-1} 
+\frac14 {E}_{21;ij}^{\bar s,\bar r; n-1 } \right. \\
&& \left. - \frac14 {E}_{21;ij}^{-\bar s,-\bar r; n-1}
\right) \frac{\bar Q_{i+\bar r,j+\bar s}^{n-1}}{Q_{i+\bar r,j+\bar s}^{\tilde
r, \tilde s; n-1} },
\end{eqnarray*}
and, finally for $|r|+|s|=3,$ we have 
\begin{eqnarray*}
A_{ij}^{rs} = \hbox{sgn }(rs) \frac{\tau}{4h^4}
\frac{
\bar Q_{ij}^{n-1} \bar Q_{i+\bar r,j+\bar s}^{n-1}}
{Q_{i+\bar r, j+\bar s}^{\tilde r,\tilde s; n-1} }
\left( {E}_{12,ij}^{\bar r,0; n-1} + {E}_{21,ij}^{0,\bar s; n-1} \right),
\end{eqnarray*}
where we have denoted $\bar r =\hbox{sgn}(r), \bar s =\hbox{sgn}(s)$ and 
$\tilde r=r-\bar r, \tilde s = s-\bar s$.


\vfill

\begin{thebibliography}{99}


\bibitem{B}
\sc M.~Bene\v s,
\em Numerical Solution for Surface Diffusion on Graphs,
\rm In: Proc. of Czech Japanese Seminar in Appl.  Math. 2005,
Bene\v s M., Kimura M. and Nakaki T., Eds.,
COE Lecture Notes, Vol. 3, Faculty of Math., Kyushu
University Fukuoka, October 2006, pp. 9--25..

\bibitem {BMN}
\sc E.~B\"ansch, P.~Morin, R.~Nochetto,
\em Surface diffusion of graphs: Variational formulation, error analysis,
and simulation,
\rm SIAM J. Numer. Anal., 42 (2004), pp.~773--799.

\bibitem {BGN}
\sc J.~W.~Barrett, H.~Garcke, R.~N\"urnberg
\em A Parametric Finite Element Method for Fourth Order Geometric Evolution
Equations
\rm Journal of Computational Physics, 222 (2007), 441--467.


\bibitem {CS}
\sc G.~Citti, A.~Sarti,
\em A cortical based model of perceptual completion in the roto-translation
space, 
\rm J. Math. Imaging and Vision, 24(3) (2006) pp.~307--326.



\bibitem {CT}
\sc J.~W.~Cahn, J.~E.~Taylor,
\em Surface motion by surface diffusion,
\rm Acta Metallica Materiala, 42 (1994),  pp.~1045--1063.

\bibitem {CDDRR}
\sc U.~Clarenz, U.~Diewald, G.~Dziuk, M.~Rumpf, R.~Rusu, 
\em A finite element method for surface restoration with
smooth boundary conditions, 
\rm Computer Aided Geometric Design, 21 (2004), pp.~427--445.

\bibitem{DG}
\sc K.~Deckelnick and H.-Ch.~Grunau, 
\em Boundary value problems for the one-dimensional Willmore equation ­ Almost
explicit solutions, \rm Preprint 2005

\bibitem{DD}
\sc K.~Deckelnick,  G.~Dziuk,
\em Error analysis of a finite element method for the Willmore flow of graphs.
\rm Interfaces Free Bound. 8, No. 1, (2006), pp.~21-46.

\bibitem {DR}
\sc M.~Droske, M.~Rumpf,
\em A level set formulation for Willmore flow,
\rm Interfaces and Free Boundaries, 6(3) (2004), pp.~361--378.

\bibitem {DKS}
\sc G.~Dziuk, E.~Kuwert, R.~Schatzle, 
\em Evolution of elastic curves in $\R^n$: existence and computation,
\rm SIAM J. Math. Anal., 33 (2002),
pp.~1228--1245.

\bibitem {E}
\sc L.~Euler, 
\em Methodus Inveniendi Lineas Curvas: Additamentum I, De Curvis Elasticis, 
\rm Opera Omnia, Z\"urich: Orell Fassli, Ser. 1, 24 (1952),  pp.~231­-297.

\bibitem{FM1}
\sc P.Frolkovi\v c, K.Mikula, 
\em Flux-based level set method: 
a finite volume method for evolving interfaces,
\rm Applied Numerical Mathematics, to appear, doi:10.1016/j.apnum.2006.06.002.

\bibitem{FM2}
\sc P.Frolkovi\v c, K.Mikula,
\em High-resolution flux-based level set method, 
\rm SIAM J. Sci. Comp., to appear.

\bibitem {Hou1}
\sc T.Y.~Hou, J.~Lowengrub, M.~Shelley,
\em Removing the stiffness from interfacial flows and surface tension,
\rm J.~Comput.~Phys., 114 (1994), pp.~312--338.

\bibitem{KWT}
\sc M.~Kass, A.~Witkin, D.~Terzopulos, 
\em Snakes: active contour models,
\rm International Journal of Computer Vision, 1 (1987), pp.~321--331.

\bibitem{K2}
\sc M.~Kimura,
\em Numerical analysis for moving boundary problems using the boundary
tracking method,   
\rm Japan J.~Indust.~Appl.~Math.,  14 (1997), pp.~373--398.

\bibitem {MS2}
\sc K.~Mikula, D.~\v Sev\v covi\v c,
\em  Evolution of plane curves driven by a nonlinear function of curvature 
and anisotropy, 
\rm SIAM J. Appl. Math., 61 (2001), pp.~1473--1501.

\bibitem {MS3}
\sc K.~Mikula, D.~\v Sev\v covi\v c,
\em A direct method for solving an anisotropic mean curvature flow of planar
curve with an external force,
\rm Mathematical Methods in Applied Sciences, 27(13) (2004) pp.~1545-1565.

\bibitem {MS_CVS}
\sc K.~Mikula, D.~\v Sev\v covi\v c,
\em  Computational and qualitative aspects of evolution of curves driven by
curvature and external force,
\rm Comput. Visual. Sci., 6 (2004), pp.~211--225.

\bibitem {MS_ALG}
\sc K.~Mikula, D.~\v Sev\v covi\v c,
\em  Tangentially stabilized Lagrangean algorithm for elastic curve evolution
driven by intrinsic Laplacian of curvature
\rm ALGORITMY 2005, Conference on Scientific Computing, 
Vysoke Tatry-Podbanske, Slovakia, March 13-18, 2005, 
Proceedings of contributed papers and posters (2005), pp.~32--41.

\bibitem{O}
\sc T.~Oberhuber,
\em Numerical Solution for the Willmore Flow of Graphs,
\rm In: Proc. of Czech Japanese Seminar in Appl.  Math. 2005,
Bene\v s M., Kimura M. and Nakaki T., Eds.,
COE Lecture Notes, Vol. 3, Faculty of Math., Kyushu
University Fukuoka, 2006, pp. 126--138.

\bibitem {Se2}
\sc J.A.~Sethian,
\em Level Set Methods and Fast Marching Methods: Evolving
Interfaces in 
Computational Geometry, Fluid Mechanics, 
Computer Vision, and Material Science,
\rm Cambridge University Press, New York, 1999.

\bibitem {SY}
\sc Y.~Saad,
\em Iterative Methods for Sparse Linear Systems (2nd edition),
\rm SIAM, 2003. 

\bibitem {ZH}
\sc H.~Zhao,
\em Fast Sweeping Method for Eikonal Equations
\rm Mathematics of Computation, 74, (2005), pp.~603--627.

\bibitem {ZC}
\sc W.~Zhu, T.~Chan, 
\em A variational model for capturing illusory contours using curvature,
\rm J. Math. Imaging and Vision, 27(1) (2007) pp.~29--40
\end{thebibliography}
\end{document}